\documentclass[12pt]{article} 

\usepackage{amsfonts}
\usepackage{amsthm}
\usepackage{xr}
\usepackage{graphicx}
\usepackage{amsmath}
\usepackage{epsfig}
\usepackage{color}

\theoremstyle{plain}

\def\d{{\partial}}

\def\C{{\mathbb C}}

\def\cC{{\cal C}}

\def\DD{{\mathbb D}}
\def\D{{\mathbb D}}
\def\cD{{\cal D}}
\def\NN{{\mathbb N}}

\def\RR{{\mathbb R}}
\def\R{{\RR}}

\def\TT{{\mathbb T}}
\def\T{{\mathbb T}}

\def\e{{\varepsilon}}
\def\cB{{\cal B}}

\def\cH{{\cal H}}
\def\trait (#1) (#2) (#3){\vrule width #1pt height #2pt depth #3pt}
\def\fin{\hfill\trait (0.1) (5) (0) \trait (5) (0.1) (0) \kern-5pt 
\trait (5) (5) (-4.9) \trait (0.1) (5) (0)}
\def\tr{\mbox{tr }}
\def\supess{\mathop{\mbox{ess sup}\,}}

\newtheorem{remarquesubsect}{Remark}[subsection]
\newtheorem{remarquesubsubsect}{Remark}[subsubsection]

\newtheorem{corollarysubsubsect}{Corollary}[subsubsection]

\newtheorem{propositionsubsect}{Proposition}[subsection]
\newtheorem{propositionsubsubsect}{Proposition}[subsubsection]

\newtheorem{lemmasubsect}{Lemma}[subsection]
\newtheorem{lemmasubsubsect}{Lemma}[subsubsection]

\newtheorem{theoremsubsect}{Theorem}[subsection]
\newtheorem{theoremsubsubsect}{Theorem}[subsubsection]

\newcommand{\Bk}{\color{black}}
\newcommand{\Rd}{\color{black}}
\newcommand{\Gr}{\color{black}}
\newcommand{\Bl}{\color{black}}

\begin{document}

\author{\begin{tabular}{cc}
Laurent Baratchart \footnote{INRIA Sophia Antipolis M\'editerran\'ee, 2004 route des
Lucioles, BP 93, 06902 Sophia-Antipolis, France. Email:
Laurent.Baratchart@sophia.inria.fr} & Juliette Leblond \footnote{INRIA
Sophia Antipolis M\'editerran\'ee. Email: Juliette.Leblond@sophia.inria.fr}\\
St\'ephane Rigat \footnote{Universit\'e de Provence, CMI, 39 rue
F. Joliot-Curie, F-13453 Marseille Cedex 13, France, and LATP, CNRS,
UMR 6632. Email: rigat@cmi.univ-mrs.fr}& Emmanuel Russ
\footnote{Universit\'e Paul C\'ezanne, Facult\'e des Sciences et
Techniques de Saint-J\'er\^ome, Avenue Escadrille Normandie-Ni\'emen,
13397 MARSEILLE Cedex 20, France, and LATP, CNRS, UMR 6632. E-mail:
emmanuel.russ@univ-cezanne.fr}
\end{tabular}}

\title{Hardy spaces of the conjugate Beltrami equation}

\maketitle

\tableofcontents

\medskip

\noindent{\small{{\bf Abstract.} We study Hardy spaces of solutions to 
the conjugate Beltrami equation with Lipschitz coefficient
on Dini-smooth simply connected
planar domains, in the range of exponents $1<p<\infty$.
We analyse their boundary behaviour and 
certain density properties of their traces.
We derive on the way an analog of the Fatou theorem
for the Dirichlet and Neumann problems associated with
the  equation $\mbox{div}(\sigma\nabla u)=0$ with 
$L^p$-boundary data.  }} 


\section{\Rd Background and motivation \Bk}
Classical Hardy spaces of the disk or the half-plane
lie at the crossroads between complex 
and Fourier analysis, and many developments in spectral theory
and harmonic analysis originate in them \cite{Nikolskiis,stein,SteinWeiss}.

From a spectral-theoretic point of view, the shift operator and its 
various compressions play a fundamental role and stress deep connections
between function theory on the one hand, control, approximation, and 
prediction theory on the other hand \cite{Nikolskii,Parfenov,Peller}.
In particular, Hankel and Toeplitz operators on Hardy classes
team up with standard functional analytic tools
to solve extremal problems where a function, given on part or all of the
boundary, is to be approximated by traces of analytic or meromorphic
functions \cite{AAK,BLPprep,bl,blp2,blp3,BarSeyf,cps,KN,Prokhorov2002}. 
Such techniques are of particular
relevance to identification and design of linear control systems
\cite{abl,blpt,DFT,Glover,KN,Peller,seyfertIMS2003}.
In recent years, on regarding the Laplace equation as a compatibility
condition for the Cauchy-Riemann system, analogous
extremal problems were set up to handle inverse Dirichlet-Neumann
issues for 2-D harmonic functions \cite{BBHL,BLMS,BMSW06,imh2,LJMP}. 
Laying grounds for a similar approach to inverse problems involving
more general diffusion elliptic equations in the plane
has been the initial motivation for the authors
to undertake the present study. The equations we have in mind are 
of the form
\begin{equation}
\label{diffus}
\mbox{div} (\sigma\nabla u)=0,~~~~\sigma~\mbox{real-valued},~~0<c<\sigma<C,
\end{equation}
which may be viewed, upon setting $\nu=(1-\sigma)/(1+\sigma)$,
as a compatibility condition for the
\emph{conjugate} Beltrami equation
\begin{equation}
\label{CBa}
\overline\partial f  = \nu  \, \overline{\partial f},~~~~
f=u+iv,~~~~\nu\in\RR,~~~~|\nu|<\kappa<1.
\end{equation}
\Gr This connection between (\ref{diffus}) and (\ref{CBa}) was instrumental in \cite{ap} for the solution of Calder\'on's conjecture in dimension $2$. \Bk Equation (\ref{CBa}) breaks up into a system of two
real equations that reduces to the Cauchy-Riemann system when $\sigma=1$.
It differs considerably from the better known Beltrami equation:
$\overline{\partial} f=\nu\,\partial f$
whose solutions are the so-called quasi-regular mappings,
which are complex linear and have been extensively studied by many authors, see
\cite{ahlfors, ast, aim, boj, cfmoz, dav, gmsv, imsurv, im2, im1, lv, mor}.

An interesting example of a free boundary inverse problem involving an equation
like (\ref{diffus}) --albeit in doubly connected geometry--
arises when trying to recover the surface of the plasma in 
plane sections of a toroidal tokamak from the so-called Grad-Shafranov equation \cite{blum}.

In the present paper, we limit ourselves to the case 
when $\sigma$ is Lipschitz-continuous. Moreover, we merely consider 
analogs $H_\nu^p$ to the classical Hardy spaces $H^p$
in the range $1<p<\infty$, on Dini-smooth simply connected 
domains.
From the perspective of harmonic analysis, the main features of Hardy space
theory in this range of exponents \cite{duren,gar}
are perhaps the Fatou theorem on 
non-tangential boundary values,
the $p$-summability of the non-tangential maximal function,
the boundedness of the conjugation operator, which is the
prototype of a convolution operator of Calder\'on-Zygmund type,
and the fact that subsets of positive measure of
the boundary are uniqueness sets 
(this is false in dimension greater than 2 \cite{BW}).
In this work, we show that Hardy solutions to (\ref{CBa}) 
enjoy similar properties,
and we use them to establish the density of the traces of 
such solutions in $L^p({\Bl \Gamma})$ whenever ${\Bl \Gamma}$ is a subset of non-full 
measure of the boundary. This fact, whose proof is straightforward
for classical Hardy spaces \cite{bl} and can be generalized
to harmonic gradients
in higher dimensions when ${\Bl \Gamma}$ is closed \cite{ABLP}, is of fundamental
importance in extremal problems with incomplete boundary data and one of
the main outcome of the paper.

The generalized Hilbert transform $\mathcal{H}_\nu$ involved in
(\ref{CBa}), that maps the boundary values of $u$ to the 
boundary values of $v$, was introduced and studied in \cite{ap,ap2} when $p=2$
for less smooth ({\it i.e.} measurable bounded) $\sigma$ but smoother
({\it i.e.} Sobolev $W^{1/2,2}$) arguments. Here, we shall prove its
$L^p$ and $W^{1,p}$ boundedness and compare it to the classical 
conjugation operator. In addition, studying its adjoint will 
lead us to a representation theorem for the dual of $H_\nu^p$
which generalizes the classical one. The latter is again
of much importance when studying extremal problems. 

On our way to the proof of the density theorem, we 
establish regularity results for the solutions of
(\ref{diffus}) which are not entirely classical. For example, 
we obtain an analog of the Fatou theorem
concerning solutions of the Dirichlet problem for (\ref{diffus}) with
$L^p$ boundary data, including $L^p$-estimates for the 
nontangential maximal function. 
Also, the gradient
of  a solution to the Neumann problem with $L^p$ data has
$L^p$ nontangential boundary values a.e. as well as
$L^p$-summable nontangential maximal function.
For the ordinary Laplacian this
is known to hold on $C^1$ domains in all dimensions
\cite{FJR}, and on Lipschitz domains for restricted range of $p$
\cite{JeKe}. But for diffusion equations of the form (\ref{diffus}),
the authors could not locate such a result in the literature, even in 
dimension two; when $p=2$ and $\sigma$ is smooth, it follows from 
\cite[{Ch. II, Thms 7.3, 8.1}]{LioMag1} that this gradient 
converges radially in $L^2$ on (parallel transportations of)
the boundary, and the result could be carried over to 
any $p\in(1,\infty)$ using the methods of \cite{Grisg},
but no pointwise estimates are 
obtained this way\footnote{In other respects
the results of \cite{LioMag1} are of course much more general
since they deal with arbitrary non-homogeneous elliptic equations
in any dimension and can handle distributional boundary conditions.}. 
  
The definition of generalized Hardy spaces
that we use --see (\ref{Smirnovdef}) below-- dwells on the existence of 
harmonic majorants
for $|f|^p$, and also on the boundedness of $L^p$ norms of $f$ on 
Jordan curves tending to the boundary of the domain \cite{duren}. 
As in the classical case, these two definitions of Hardy spaces 
coincide on Dini-smooth domains (the only case of study below)
but not over non-smooth domains --where arclength on the boundary and
 harmonic measure are no longer mutually absolutely continuous.
 
Although equation (\ref{CBa}) is real linear only, our methods
of investigation rely on complex analytic tools. In particular, we elaborate
on ideas and techniques from \cite{bn} and we use standard facts from
classical Hardy space theory together with well-known
properties of 
the Beurling transform. This entails that higher dimensional analogs 
of our results, if true at all, require new ideas to be proven.

We made no attempt at expounding the limiting cases $p=1,\infty$.
These have generated the deepest developments in the classical theory,
centering around BMO and Fefferman duality,
but trying to generalize them would have made the paper unbalanced and they
are left here for further research.

Finally, we did not consider Hardy spaces over
doubly connected domains, in spite of the fact 
that the above-mentioned application to free boundary problems in 
plasma control takes place in an annular geometry. 
Including these would have made for a lengthy paper,
but the results below lay ground for such a study.

\section{\Rd Notations for function spaces \Bk}
\label{sec:pb}
Throughout, $\D$ is the open unit disk and 
$\T$ the unit circle in the complex plane 
$\C$. We let $D_r$ and $\T_r$ stand for 
the open disk and the circle centered at $0$  with radius $r$.
For $I$ an open subset of $\T$,
endowed with its natural differentiable 
structure, we put $\cD(I)$ for the space of
$C^\infty$ complex functions supported on $I$.\par

\medskip

\noindent If $\Omega\subset \C$ is a smooth domain (the meaning of ``smooth'' will be clear from the context), we say that a sequence $\xi_n\in\Omega$ approaches
$\xi\in\partial\Omega$ non tangentially
if it converges to $\xi$ while no limit point of $(\xi_n-\xi)/|\xi_n-\xi|$
is tangent to $\partial\Omega$ at $\xi$.
A function $f$ on $\Omega$
has non tangential 
limit $\ell$ at $\xi$ if $f(\xi_n)$ tends to $\ell$ for any sequence 
$\xi_n$ which approaches  $\xi$ non tangentially. 

\subsection{\Rd H\"older spaces \Bk}
If $\Omega\subset\RR^2$ is open, 
$C^{k,\gamma}(\overline{\Omega})$ indicates the subspace of complex
functions whose derivatives are bounded and continuous up to order $k$, 
while those
of order $k$ satisfy a H\"older condition
of exponent $\gamma\in(0,1]$. Such functions extend continuously to 
$\overline{\Omega}$ together with their derivatives of order at most $k$.
A complete
norm on $C^{k,\gamma}(\overline{\Omega})$
is obtained by putting
\[\left\Vert f\right\Vert_{C^{k,\gamma}(\overline{\Omega})}:=
\sup_{0\leq|\lambda|\leq k}\|f^{(\lambda)}\|_{L^\infty(\Omega)}+
\sup_{\stackrel{|\lambda|= k}{\xi\neq\zeta}}
\frac{|f^{(\lambda)}(\xi)-f^{(\lambda)}(\zeta)|}{|\xi-\zeta|^\gamma},
\]
where $\lambda=(\lambda_1,\lambda_2)$ is a multi-index, $|\lambda|=\lambda_1+\lambda_2$,
and $f^{(\lambda)}$ is the corresponding derivative.
The space $C^\infty(\overline{\Omega})=\cap_k
C^{k,1}(\overline{\Omega})$ of smooth functions 
up to the boundary of $\Omega$ is topologized with the seminorms
$\|~\|_{C^{k,1}(\overline{\Omega})}$, where $k$ ranges over $\NN$.
We put $C_{loc}^{k,\gamma}(\Omega)$
for the functions whose restriction to any relatively compact open subset
$\Omega_0$ of $\Omega$ lies in $C^{k,\gamma}(\overline{\Omega_0})$. 
A family  of semi-norms making $C_{loc}^{k,\gamma}(\Omega)$ 
into a Fr\'echet space is given by $\|~\|_{C^{k,\gamma}(\Omega_n)}$, 
with $\Omega_n$ a sequence of relatively compact open subsets 
exhausting $\Omega$.

\subsection{\Rd Lebesgue and Sobolev spaces \Bk}
We coordinatize $\R^2 \simeq \C$ by $\xi=x+iy$
and denote interchangeably the (differential of) planar
Lebesgue  measure by 
\[dm (\xi) = dx \, dy = (i/2) d\xi\wedge d\overline{\xi}\,,\]
where $d\xi=dx+i dy$ and $d\overline{\xi}=dx-i dy$.
For $1\leq p\leq +\infty$ and $E$ a measurable subset of $\C$, 
we put $L^p(E)$ for the familiar Lebesgue space with respect to $dm$.

\noindent Let $\mathcal D(\Omega)$ 
be the space of complex $C^\infty$ functions with compact support in 
$\Omega$,  endowed with the inductive topology. Its dual 
$ {\mathcal D}^{\prime}(\Omega)$ is the usual space 
of distributions on $\Omega$.
Whenever $T\in {\mathcal D}^{\prime}(\Omega)$, 
we use the standard notations:
\[
\partial T=\partial_zT=\frac 12\left(\partial_x-i\partial_y\right)T
~~\mbox{ and }~~\overline{\partial}T=\partial_{\overline{z}}T=\frac 12\left(\partial_x+i\partial_y\right)T,
\]
and record the obvious identity: $\overline{\d T}=\bar\d\,\overline{T}$.

We denote by $W^{1,p}(\Omega)$ the Sobolev space of 
those $f\in L^p(\Omega)$ whose distributional derivatives $\partial f $ and
$\overline{\partial}f$ also belong to $L^p(\Omega)$. The norm on 
$W^{1,p}(\Omega)$ is defined by
\[\|f\|^p_{W^{1,p}(\Omega)}:=\|f\|_{L^p(\Omega)}^p+\|\partial f\|_{L^p(\Omega)}^p+\|\overline{\partial} f\|_{L^p(\Omega)}^p.\]
\Rd When $\Omega$ is smooth, any function $f\in W^{1,p}(\Omega)$ has a trace on $\partial\Omega$ (all the domains under consideration in the present paper are smooth enough for this to be true) which will be denoted by ${\mbox tr }f$. \Bk\par
\noindent The symbol $W^{1,p}_{loc}(\Omega)$ refers to those distributions whose 
restriction to any relatively compact open subset
$\Omega_0$ of $\Omega$ lies in $W^{1,p}(\Omega_0)$. 
Equipped with the semi-norms $\|~\|_{W^{1,p}(\Omega_n)}$,
$W^{1,p}_{loc}(\Omega)$ is a Fr\'echet space.

\par
\noindent The space $W^{1,\infty}(\Omega)$ is isomorphic to 
$C^{0,1}(\overline{\Omega})$ 
\cite[Ch.VI, Sec. 6.2]{stein}. In particular, every $f\in W^{1,\infty}(\Omega)$ 
extends (Lipschitz) continuously to the boundary
$\partial\Omega$ of $\Omega$.
\par
\noindent A $C^1$-smooth Jordan curve is the injective image of $\T$ under a
nonsingular continuously differentiable map from $\T$ into $\C$. 
For $\Gamma$ an open subset of such a curve,
we denote by $L^p(\Gamma)$ the Lebesgue space
with respect to (normalized) arclength (\Rd there should be no confusion with our previous notation
$L^p(E)$, as the context will always remain clear\Bk) and by
$W^{1,p}(\Gamma)$ the Sobolev space of those absolutely continuous
$\varphi\in L^p(\Gamma)$ whose tangential derivative $\partial_t \varphi$
with respect to arclength again lies in $L^p(\Gamma)$. 
A complete norm on $W^{1,p}(\Gamma)$ is 
given by 
\[\|\varphi\|^p_{W^{1,p}(\Gamma)}:=\|\varphi\|_{L^p(\Gamma)}^p+
\|\partial_t \varphi\|_{L^p(\Gamma)}^p.
\]
Note that, on $\T$,  the tangential derivative $\partial_t h$ 
coincides with $\partial_\theta h/(2\pi)$ where
$\partial_\theta h$ indicates the derivative with respect to $\theta$
when $\varphi$ is written as a function of $e^{i\theta}$. 
\par
\noindent

We shall have an occasion to deal with $W^{2,p}(\Omega)$, 
comprised of $W^{1,p}(\Omega)$-functions whose first derivatives again lie in 
$W^{1,p}(\Omega)$. A norm on $W^{2,p}(\Omega)$ is obtained by setting 
\[\|f\|^p_{W^{2,p}(\Omega)}=\|f\|^p_{L^p(\Omega)}+\|\partial f\|_{W^{1,p}(\Omega)}^p+\|\overline{\partial} f\|_{W^{1,p}(\Omega)}^p.\]

As is customary, we indicate with a subscript ``0'' the
closure of $C^\infty$ compactly supported  functions in an ambient space.

\par

\medskip

\noindent \Rd Finally, we indicate with a subscript ``$\R$'', like in
$L^p_{\R}({\Bl\Omega})$,
the real subspace of real-valued functions in a given space. 

\section{\Rd Definition of Hardy spaces \Bk}
\subsection{\Rd An elliptic equation \Bk}

In the present paper, we investigate the $L^p$ boundary behaviour of solutions 
to a second order elliptic equation in divergence form on a planar domain.
More precisely, \Rd let $\Omega\subset \R^2$ be a smooth simply connected domain (most of the time, we will take $\Omega=\D$, except in Section \ref{dini-sec}, where $\Omega$ will be assumed to be Dini-smooth) \Bk
and $\sigma\in W^{1,\infty}(\Omega)$ be such that, for two constants 
$c,C$,  one has
\begin{equation} \label{ellipticsigma}
0<c \leq \sigma \leq C\,.
\end{equation}
With the standard notation $\nabla u:=(\partial_x u\,,\,\partial_y u)^T$ and 
$\mbox{div }g=\partial_xg+\partial_y g$, where the superscript ``$T$'' 
means ``transpose'', the elliptic equation that we will consider is
\begin{equation} \label{div}
\mbox{div} (\sigma\nabla u)=0 \ \ \mbox{ a.e. in }\Omega,
\end{equation}

\noindent Our approach to (\ref{div}) proceeds {\it via} the study
of a \emph{complex} elliptic equation of the first order, namely
the {\it conjugate} Beltrami equation:
\begin{equation}
\label{dbar} 
\overline\partial f  = \nu  \, \overline{\partial f} \ \ \mbox{ a.e. in } \Omega \, ,
\end{equation}
where $\nu\in W^{1,\infty}(\Omega)$ is a {\it real valued} function
that satisfies:
\begin{equation} \label{kappa} 
\left\Vert \nu\right\Vert_{L^\infty(\Omega)}\leq \kappa \ \ \mbox{ for some }
\kappa\in (0,1) \, .
\end{equation}
Formally, equation (\ref{dbar}) decomposes into a system of two real
elliptic equations of the second order in divergence form. Indeed,
for $f=u+iv$ a solution to (\ref{dbar}) with real-valued $u$, $v$,
we see on putting $\sigma=(1-\nu)/(1+\nu)$ that
$u$ satisfies equation (\ref{div})
while $v$ satisfies
\begin{equation} \label{system2}
\mbox{div}\left(\frac 1{\sigma}\nabla v\right)=0\ \ \mbox{a.e. in }\Omega.
\end{equation}
Note also, from the definition of $\sigma$,  that (\ref{kappa}) implies
(\ref{ellipticsigma}).
Conversely, let 
$u$ be a real-valued solution to (\ref{div}). Then, since 
$\partial_y(-\sigma\partial_yu) = \partial_x(\sigma\partial_xu)$ 
and $\Omega$ is simply connected, there is a real-valued function $v$,
such that
\begin{equation} \label{system}
\left\{
\begin{array}{l}
\partial_xv=-\sigma\partial_y u,\\
\partial_yv=\sigma\partial_xu,
\end{array}
\right.
\end{equation}
hence  $f=u+iv$ satisfies (\ref{dbar}) with 
$\nu=(1-\sigma)/(1+\sigma)$. Moreover, 
(\ref{ellipticsigma}) implies (\ref{kappa}).

In the present work, we consider several classes of solutions 
to (\ref{dbar}) for which the formal manipulations
above  will be given a precise meaning. All classes we shall deal with are
embedded in $L^p(\Omega)$ for some $p\in(1,\infty)$, in which case
the solutions to \eqref{div}, (\ref{dbar}), and \eqref{system2}
can be understood in the distributional sense. This only requires defining
distributions like $\sigma\partial_x u$, which is done 
naturally using Leibniz's rule\footnote{For instance if $f \in L^p(\Omega)$ and 
 $\nu\in W_\R^{1,\infty}({\Bl \Omega})$,  we  define by  
$\nu \overline{\partial  f}$ to be the distribution 
 \[
 \langle \nu \overline{\partial  f}, \phi \rangle = - \int_\Omega 
(\nu \overline{f}\, \overline{\partial} \phi  + \overline{\partial} \nu 
 \overline{f} \phi) dm \, , \quad\forall \phi \in {\cal D}(\Omega) \,. 
 \] 
   }
when $\sigma\in W^{1,\infty}(\Omega)$ and 
$u\in L^p(\Omega)$
\cite{cfmoz}. 

It will turn out that our solutions actually lie  
in $W^{1,p}_{loc}(\Omega)$ for some $p\in(1,\infty)$, in which case 
\eqref{dbar} may as well be interpreted in the pointwise  sense while
\eqref{div} becomes  equivalent to
\begin{equation}
\label{weaksol}
\int_\Omega \sigma\nabla u.\nabla g\,dm =0,\ \ \ \ g\in\cD_\R(\Omega),
\end{equation}
where the dot indicates the Euclidean scalar product in $\R^2$. 
This follows easily from the fact that the product of a function in $W^{1,\infty}(\Omega)$ by a function in 
$W^{1,p}_{loc}(\Omega)$ again lies in $W^{1,p}_{loc}(\Omega)$ and its distributional derivative can be computed according to the Leibniz rule.
We shall make use of these observations without further notice.\par

\medskip

\noindent To find $u$ with prescribed trace on $\partial\Omega$ is known as the 
Dirichlet problem for \eqref{div} in $\Omega$.  
In light of the previous discussion, we slightly 
abuse terminology and still refer to the issue of finding $f$
with prescribed 
$\mbox{\rm Re}f$ on $\partial\Omega$ as being
the Dirichlet problem for \eqref{dbar}.\par

\medskip

\noindent For simplicity, we shall work entirely over the unit disk $\D$
and only later, in Section \ref{dini-sec}, shall we 
indicate how one can carry our results over to Dini-smooth domains.
As became customary in analysis, we tend to use the same symbol to 
mean possibly 
different constants, with subscripts indicating the dependence of the constant
under examination.

\par
\noindent

%
%

\bigskip

\par
When $\nu\in W^{1,\infty}(\D)$, the solvability in $W^{1,p}(\D)$
of the Dirichlet problem
for \eqref{dbar} with boundary data in \Rd the fractional Sobolev space $W^{1-1/p,p}_{\R}(\T)$ (an intrinsic definition of which 
can be found in \cite[{ Thm} 7.48]{Adams}) \Bk is a straightforward consequence of the known 
solvability of the corresponding Dirichlet problem for equation (\ref{div})
\cite{campanato}. 
\Rd We shall however state and establish this fact which is 
our point of departure (see Theorem \ref{dirichlet} in Section \ref{statement} below). \Bk\par

\medskip

\noindent Below, we relax the assumptions on the boundary data, assuming only
they belong to $L^p(\T)$. Of course, the solution of the Dirichlet problem 
will no longer belong to $W^{1,p}(\D)$ in general, but rather 
to  some generalized  
Hardy space $H^p_\nu(\D)$ that we shall define and study throughout the paper,
starting in the next section.

\subsection{Definition of Hardy spaces}
\label{notation}
For $1<p<\infty$, we denote by  $H^p(\D)$ the  classical Hardy 
space of holomorphic functions $f$ on $\D$ such that 
\begin{equation}
\label{Smirnovdef}
\left\Vert f\right\Vert_{H^p(\D)}:=\supess_{0<r<1}  \left\Vert f
\right\Vert_{L^p(\T_r)} <+\infty \,, 
\end{equation}
where
\[
\left\Vert f\right\Vert_{L^p(\T_r)} := 
\left(\frac
1{2\pi}\int_0^{2\pi} \left\vert
f(re^{i\theta})\right\vert^pd\theta\right)^{1/p} \, ,
\]
 by our convention that arclength gets normalized, see \cite{duren, gar}.

Of course $H^p$ can be introduced for $p=1,\infty$ as well,
but we do not consider such exponents here.
We extend the previous definition to
two classes of generalized analytic functions as follows.
\par
\subsubsection{The class $H^p_{\nu}(\D)$}
\noindent If $\nu\in W^{1,\infty}_{\R}(\D)$ satisfies (\ref{kappa}),
and $1<p<+\infty$, 
we define a generalized Hardy space $H_\nu^p(\D)$
to consist of those Lebesgue measurable functions $f$ on $\D$  
satisfying 
 \begin{equation} \label{esssuplp}
 \left\Vert f\right\Vert_{H^p_{\nu}(\D)}:= \supess_{0<r<1}
 \left\Vert f\right\Vert_{L^p(\T_r)} <+\infty \, 
 \end{equation}
 that solve \eqref{dbar} 
 in the sense of distributions  on $\D$; note that \eqref{esssuplp}
implies $f\in L^p(\D)$.
It is not difficult to see that $\|.\|_{H^p_\nu(\D)}$ is a norm
making $H^p_{\nu}(\D)$ into a real Banach space. 

When $\nu=0$, then $H^p_{\nu}(\D)=H^p(\D)$ viewed as a 
\emph{real} vector space.

\par
 
\medskip As we will see in Proposition \ref{cor:propHpnu}, each 
$f\in H^p_{\nu}(\D)$ has a non-tangential limit 
a e. on $\T$ that we call the trace of $f$,
denoted by $\tr\, f$ (see Section \ref{sec:pb} for the definition of the non-tangential limit). \Rd This definition causes no discrepancy since, as we shall see in Proposition \ref{hardysobol} below, any solution of \eqref{dbar} in $W^{1,p}(\D)$ belongs to $H^p_{\nu}(\D)$ and, for an arbitrary function $f\in W^{1,p}(\D)$, the nontangential limit of $f$, when it exists, coincides with the trace of $f$ in the Sobolev sense. \Bk
It turns out that, for all $f\in H^p_{\nu}(\D)$, $\tr\, f$ lies in $L^p(\T)$
and $\|\tr\, f\|_{L^p(\T)}$ defines an equivalent norm on
$H^p_\nu(\D)$. \Bk

\medskip
\par
\noindent We single out the subspace $H^{p,0}_{\nu}$ of $H^p_{\nu}$  
consisting of those $f$ for which
\begin{equation}
\label{normal:f}
\int_0^{2\pi}
\mbox{Im} \, (\tr\, f (e^{i \theta})) \, d  \theta = 0 \, 
\end{equation}
holds.
We further let $H_{\nu}^{p,00}$ be the subspace of 
functions $f\in H_{\nu}^{p,0}$  such that
\begin{equation}
\label{normal:fv}
\int_0^{2\pi}
\tr\, f (e^{i \theta}) \, d  \theta = 0 \, .
\end{equation}
\par
\begin{remarquesubsect}
In what follows, we make use of both $H^p_{\nu}(\D)$ and $H^p(\D)$.
For simplicity, we drop the dependence on $\D$ and
denote them by $H^p_{\nu}$ and $H^p$, respectively.
In particular, $H^p$ (no subscript) 
always stands for the classical holomorphic Hardy space of the disk.
\end{remarquesubsect}
\subsubsection{The class $G^p_{\alpha}(\D)$}
For $\alpha \in L^{\infty}(\D)$ and $1<p<\infty$ (note that $\alpha$ may be complex-valued here), 
we  define another space 
$G_\alpha^p(\D)= G_\alpha^p$, consisting of those Lebesgue measurable 
functions $w$ on $\D$ such that:
$$ 
\left\Vert w\right\Vert_{G_\alpha^p}:= 
\supess_{0<r<1}
\left\Vert w \right\Vert_{L^p(\T_r)} <+\infty 
$$
and
\begin{equation}
\label{eq:w}
{\overline{\partial}} w = \alpha \overline{w} \, 
\end{equation}
in the sense of distributions on $\D$. 
Note that $\|.\|_{H^p_\nu(\D)}$ and $\|.\|_{G^p_\alpha(\D)}$ formally 
coincide, 
but the equations \eqref{dbar} and \eqref{eq:w} are different.
Again $\|.\|_{G^p_\alpha(\D)}$ makes
$G^p_{\alpha}$ into a real Banach space. 
The reason why we introduce $G^p_\alpha$ is the tight connection it has
with $H_{\nu}^p$ when we set $\alpha:=-\bar\partial\nu/(1-\nu^2)$,
as shown in Proposition \ref{trick-hardy} below. From equation 
\eqref{alphabeta} below,  we
see that $\alpha$ has this form for some $\nu\in W_\R^{1,\infty}$ 
meeting \eqref{kappa} if, and only if
$\alpha=\bar\partial h$ for some $h\in W_\R^{1,\infty}(\D)$. Making
such an assumption in the definition 
would be artificial, since most of the properties
of $G^p_\alpha$ to come are valid as soon as $\alpha\in L^\infty(\D)$. However,
if \eqref{alphabeta} holds \emph{and only in this case} (see section \ref{tricksec} below),
we shall find it convenient to introduce the space
$G^{p,0}_{\alpha}$ of those $w \in G^p_\alpha$ normalized by 
\begin{equation}
\label{normal:w}
\frac{1}{2\pi}\int_0^{2\pi} \left(\sigma^{1/2} \mbox{Im} \, {\tr\, w }\right)(e^{i \theta}) \, d   
\theta = 0 .
\end{equation}

\subsubsection{The link between $H^p_{\nu}$ and $G^p_{\alpha}$} 
\label{tricksec}
The explicit connection between $H^p_{\nu}$ and $G^p_{\alpha}$ 
is given by the following result, which relies on a transformation 
introduced in  \cite{bn}:
\begin{propositionsubsubsect} \label{trick-hardy}
Let $\nu\in W_\R^{1,\infty}(\D)$ satisfy \eqref{kappa} and define
$\sigma\in W_\R^{1,\infty}(\D)$, $\alpha\in 
L^\infty(\D)$ by 
\begin{equation} \label{alphabeta}
\sigma=\frac{1-\nu}{1+\nu},~~~~~~\alpha:=-\frac{\bar\partial\nu}{1-\nu^2} =\frac{\bar\partial\sigma}{2 \sigma} = \bar\partial \log \sigma^{1/2}.
\end{equation}
Then $f\in L^p(\D)$ solves \eqref{dbar} in the distributional sense if,
and only if \Rd the function $w$, defined by 
\begin{equation}
\label{correspw}
w :=( f - \nu \overline{f})/\sqrt{1-\nu^2} 
= \sigma^{1/2}\, u+ i \, \sigma^{-1/2}\, v \Bk
\end{equation}
does for \eqref{eq:w}. Moreover,
\begin{itemize}
\item[$(a)$]
$f  = u+ i \, v $ lies in $H_{\nu}^p$ (resp. $H^{p,0}_{\nu}$) if, and only if
the function \Rd $w$ given by (\ref{correspw}) \Bk
lies in $G_\alpha^p$ (resp. $G^{p,0}_{\alpha}$).
 \item[$(b)$]
$f\in W^{1,p}(\D)$ solves \eqref{dbar} if, and only if $w $ given by
\eqref{correspw} solves \eqref{eq:w} in $ W^{1,p}(\D)$.
\end{itemize}
\end{propositionsubsubsect}
The proof is a straightforward computation,
using that the distributional derivatives of
$(f-\nu\overline{f})/\sqrt{1-\nu^2}$ can be computed by Leibniz's rule
under our  standing assumptions, \Rd and the fact that (\ref{correspw}) can also be rewritten as
\begin{equation}
\label{correspf}
\Rd f =\frac{ w + \nu \overline{w}}{\sqrt{1-\nu^2}}. \Bk
\end{equation}
\Bk Observe that every constant $c \in \C$ is a solution to \eqref{dbar}, the
associated  $w$ {\it via} \eqref{correspw} being 
$\sigma^{1/2}\, \mbox{Re } \, c + i \, \sigma^{-1/2} \, \mbox{Im } \, c$, 
which lies in $W^{1,\infty}(\D)$ and solves (\ref{eq:w}).

\section{\Gr Statement of the results. \Bk} \label{statement}

\emph{Throughout, we assume $1 < p < \infty$, and 
we let $\nu \in W^{1, \infty}_\R(\D)$ satisfy $\left\Vert
  \nu\right\Vert_{L^\infty(\D)}\leq \kappa<1$.}
\vskip.5cm

\subsection{\Rd Solvability in Sobolev spaces \Bk}
\Rd Our first result deals with the solvability of the Dirichlet problem for \eqref{dbar} with boundary data in $W^{1-1/p,p}(\T)$: \Bk
\begin{theoremsubsect}
\label{dirichlet}
Let $p\in(1,+\infty)$ and $\nu\in W^{1,\infty}_\R(\D)$ satisfy \eqref{kappa}.
\begin{itemize}
\item[$(a)$]
To each $\varphi \in W^{1-1/p,p}_\RR(\T 
)$, there is  $f
\in W^{1,p}(\D)$ solving  \eqref{dbar}  in $\D$ and such that 
$\mbox{\rm Re} \left(\mbox{tr }f\right) = \varphi$ on $
\T$. Such an $f$ is unique up to an additive pure imaginary constant.
\item[$(b)$]
There exists $C_{p,\nu}>0$ such that the function $f$ in $(a)$,
when normalized by \eqref{normal:f},
satisfies
\begin{equation} \label{estimf}
\left\Vert f\right\Vert_{W^{1,p}(\D)}\leq C_{p,\nu} \left\Vert \varphi\right\Vert_{W^{1-1/p,p}(\T)}.
\end{equation}
\end{itemize}
\end{theoremsubsect}
\begin{remarquesubsect}
Although we will not use it,
let us point out that Theorem
\ref{dirichlet} still holds if we merely assume $\nu\in VMO(\D)$,
provided \eqref{dbar} is understood in the pointwise sense. The proof
is similar, appealing to \cite{aq} rather than \cite{campanato} to solve the
Dirichlet problem for \eqref{div}.
\end{remarquesubsect}
\begin{remarquesubsect}
\label{remNeuman}
When $\varphi\in W_\R^{1-1/p,p}(\T)$ and
$u\in W^{1,p}_\R(\D)$ is the solution to \eqref{weaksol} such that
$\tr\, u=\varphi$ granted by \cite{campanato}, 
the normal derivative $\partial_n u$
is classically \emph{defined} as the unique member of {\Bl the dual space} $W^{-1/p,p}_\R(\T)=\left(W^{1-1/q,q}_\R(\T)\right)^*$ such that 
\begin{equation}
\label{dern}
\langle \partial_n u\,,\,\sigma\psi\rangle\,=
\,\int_{\D} 
\sigma\nabla u.\nabla g\,dm  \,
,\ \ \ \ \psi\in W^{1-1/q,q}_\R(\T),\  g \in W^{1,q}(\D),\ \tr g=\psi,
\end{equation}
where $g$ is any representative of the coset $\tr^{-1} \psi$ in
$W^{1,p}_\R(\D)/W^{1,p}_{0,\R}(\D)$.  That $\partial_n u$ is well-defined {\it via} \eqref{dern}
depends on the fact that $M_\sigma$, the multiplication by $\sigma$, is an isomorphism of $ W^{1-1/q,q}_\R(\T)$;
this follows by interpolation since $M_\sigma$ is an isomorphism both of $L^q_\R(\T)$ and $W^{1,q}_\R(\T)$.
Now, if $f=u+iv
\in W^{1,p}(\D)$ is a solution to \eqref{dbar}  such that
$\mbox{\rm Re} \left(\mbox{tr }f\right) = \varphi$ as provided by Theorem \ref{dirichlet},
it is a straightforward consequence of \eqref{system} that $\partial_n u=(\partial_\theta \tr v)/\sigma$.
\end{remarquesubsect}

The results below generalize to $H_\nu^p$ and $G^p_{\alpha}$,
defined in Section \ref{notation},
some fundamental properties of holomorphic Hardy classes 
\cite{duren, gar}. Observe that, on $\D$, as the above  definition shows (see Section \ref{sec:pb}),
$f$ has a non tangential (``n.t.'') limit
$\ell$  at $e^{i\theta}\in\T$ if, and only if for every
$0<\beta<\pi/2$, 
$f(z)$ tends to $\ell$ as
$z\to e^{i\theta}$ inside any sector $\Gamma_{e^{i\theta},\beta}$
with vertex $e^{i\theta}$, of angle
$2\beta$, which is symmetric with respect to the ray $(0,e^{i\theta})$. 
The non-tangential maximal function of $f$ at $\xi\in\T$ 
is
\begin{equation}
\label{defntmax}
{\mathcal M}_f(\xi):=\sup_{z\in \D\cap\Gamma_{\xi,\beta}}
|f(z)|,
\end{equation}
where we dropped the dependence of ${\mathcal M}_f$ on $\beta$.
\par
\noindent We first mention properties of the class $G^p_{\alpha}$, 
from which those of the class $H^p_{\nu}$ will be deduced using
Proposition \ref{trick-hardy}.

\subsection{Properties of $G^p_{\alpha}$} \label{gpalphaprop}
We fix $\alpha\in L^\infty(\D)$. 
\noindent To proceed with the statements, we need to introduce
two operators that will be of constant use in the paper. 
First, for $ \psi \in L^1(\T)$, we define a holomorphic function
in $\D$ through the Cauchy operator:
\[
{\cal C} \psi(z)=\frac1{2\pi i} \int_{\T} \frac{\psi(\xi)}{\xi - z} d
\xi \, , \ z \in \D \, .
\]
It follows from a theorem of M. Riesz that $\cC$ maps $L^p(\T)$ onto
$H^p$, see the discussion after \cite[Ch. 3, { Thm} 1.5]{gar}; 
this would fail if $p=1,\infty$.

Second, for $p\in(1,+\infty)$ and
$w\in L^p(\D)$, we define
\[
Tw(z)=\frac1{2\pi i}\iint_{\D} \frac{w(\xi)}{\xi-z}d\xi\wedge d\overline{\xi}
 \, , \ z \in \D \, .
\]
The following 
representation theorem for functions in $G^p_{\alpha}$
was implicit in \cite{bn} for continuous $W^{1,2}(\D)$-solutions to 
\eqref{eq:w}:
\begin{theoremsubsect}
\label{trick2}
Let $w \in L^p(\D)$ be a distributional solution to 
\eqref{eq:w}. Then $w$ can be represented as
\begin{equation}
\label{factBN}
w(z) = \exp (s(z)) \, F(z)  \, , \ z \in \D \, ,
\end{equation}
where $s \in W^{1,l}(\D)$ for all $l\in (1,+\infty)$  and  $F$ is 
holomorphic in $\D$. Moreover, $s$ can be chosen such that its real part 
(or else its imaginary part) is 0 on $\T$ and 
\begin{equation} \label{estims} 
\|s\|_{L^\infty(\D)} \leq 4\|\alpha\|_{L^\infty(\D)}.
\end{equation}
In particular $w\in W^{1,l}_{loc}(\D)$ for all $l\in (1,+\infty)$,
and $w\in G^p_\alpha$ if, and only if
$F\in H^p$ in some, hence any factorization of the form \eqref{factBN}. 
Moreover, $w \in L^{p_1}(\D)$, for all $p_1 \in [p, 2p)$.
\end{theoremsubsect}
\begin{remarquesubsect} 
\label{rems}
By the Sobolev imbedding theorem (\Gr \cite[{ Thm 5.4, Part II}]{Adams}), \Bk
$s\in C^{0,\gamma}(\overline{\D})$ and $w\in C^{0,\gamma}_{loc}(\D)$
for all $\gamma\in (0,1)$. 
\end{remarquesubsect}
Theorem \ref{trick2} will allow for us to
carry over to $G_\alpha^p$ the essentials of the boundary behaviour of
holomorphic Hardy functions:
\begin{propositionsubsect} \label{gc-hardy} 
\par
\begin{itemize}
\item[1.]
If $w\in G^p_{\alpha}$,  then
\begin{equation}
\label{limntG}
\tr\, w(e^{i \theta}) :=\lim_{\xi\in \D,\ \xi\rightarrow e^{i \theta} \mbox{
    n.t.}} w(\xi)
\end{equation}
exists for almost every $\theta$ and
\begin{equation}
\label{Fatou}
\left\Vert \tr\, w\right\Vert_{L^{p}(\T)} \leq 
\left\Vert w\right\Vert_{G_\alpha^p} \leq c_{\alpha} \, \left\Vert \tr\,
  w\right\Vert_{L^{p}(\T)}\,  
\end{equation}
for some $c_{\alpha} >0$. Moreover,
\begin{equation} \label{limgp}
\lim_{r\rightarrow 1} \int_0^{2\pi} \left\vert w(re^{i\theta})-\mbox{tr }w(e^{i\theta})\right\vert^pd\theta=0
\end{equation}
and, for any aperture $\beta \in (0, \pi/2)$ of the sectors 
used in definition
 \eqref{defntmax}, there is a constant $C_{p,\alpha,\beta}$ such that
\begin{equation}
\label{maxinegw}
\|{\mathcal M}_w\|_{L^p(\T)}\leq C_{p,\alpha,\beta}\|\tr\, w\|_{L^p(\T)}.
\end{equation}
\item[2.] 
If $w\in G_{\alpha}^p$ and $w\not\equiv0$, then 
$\log |\tr\, w| \in L^1(\T)$; moreover the zeros of $w$ are isolated 
in $\DD$, and  numbering them as $\alpha_1,\alpha_2,\cdots$, counting repeated 
multiplicities,  it holds that
\begin{equation}
\label{blaschke}
\sum_{j=1}^\infty (1-|\alpha_j|)<+\infty.
\end{equation}
\item[3.]
Let $w\in L^p(\D)$. Then $w\in G_\alpha^p$ if,
and only if  there is a function $\varphi\in L^p(\T)$ such that
\begin{equation}
\label{wCT}
w =  \cC \varphi + T(\alpha\overline{w}) \, , \mbox{a.e. in } \D \, .
\end{equation}
In this situation,
\begin{equation} \label{gpnorm}
\left\Vert w\right\Vert_{G^p_{\alpha}}\leq C_{p,\alpha}\left(\left\Vert w\right\Vert_{L^p(\D)}+\left\Vert \varphi\right\Vert_{L^p(\T)}\right).
\end{equation}
A valid choice in \eqref{wCT} is 
$\varphi=\tr\, w$.
\item[4.]
If $w\in G^{p}_{\alpha}$  satisfies \eqref{normal:f}  and 
$\mbox{Re tr }w=0$ a.e. on $\T$, 
then $w\equiv 0$ in $\D$. When \eqref{alphabeta} holds, the same 
is true if $w\in G^{p,0}_{\alpha}$. 
\end{itemize}
\end{propositionsubsect}
\begin{remarquesubsect}
\label{rmkY}
From \eqref{Fatou} and the completeness of $G_\alpha^p$, we see 
that $\tr\, G_\alpha^p$
is a closed subspace of $L^p(\T)$.
We also observe, in view of the M. Riesz theorem,
that assertion 3 can be recaped as:
$w \in G_\alpha^p
\Longleftrightarrow w - T(\alpha \overline{w})  \in H^p$.
\end{remarquesubsect}
Theorem \ref{trick2} and Proposition \ref{gc-hardy} 
are proven  in Section \ref{secprop:w}.

\subsection{Properties of the Hardy class $H_{\nu}^p$} \label{kkkkk}

\begin{propositionsubsect}
\label{cor:propHpnu}
The following statements hold true.
\begin{itemize}
\item[(a)]  If $f\in H_{\nu}^p$, then $f$ has a non-tangential limit 
{\it a.e.} on  $\T$, denoted by $\mbox{tr }f$,
the $L^p(\T)$-norm of which is equivalent to the $H^p_{\nu}$-norm of $f$:
\begin{equation}
\label{Fatouf}
\left\Vert \tr\, f\right\Vert_{L^{p}(\T)} \leq 
\left\Vert f\right\Vert_{H^p_\nu(\D)} \leq c_\nu\, \left\Vert \tr\,
  f\right\Vert_{L^{p}(\T)}\, .
\end{equation}
Moreover
\[
\lim_{r\rightarrow 1} \int_0^{2\pi} \left\vert f(re^{i\theta})-\mbox{tr }f(e^{i\theta})\right\vert^pd\theta=0,
\]
and we have that $f\in L^{p_1}(\D)$ for $p\leq p_1<2p$.
\item[(b)] The image space $\tr\, H_{\nu}^p$ (resp. $\tr\, H_{\nu}^{p,0}$)
is closed in $L^p(\T)$.
\item[(c)] Each $f\in H_{\nu}^p$ is such that  $\log |\tr\, f| \in L^1(\T)$ 
unless $f \equiv 0$.
\item[(d)] If $f\in H_{\nu}^p$ and $f\not\equiv0$, then its zeros are isolated 
in $\DD$; if we enumerate them as $\alpha_1,\alpha_2,\cdots$, 
counting repeated multiplicities, then \eqref{blaschke} holds.
\item[(e)] For any aperture $\beta \in (0, \pi/2)$ of the sectors used in definition
\eqref{defntmax}, there is a constant $C_{p,\nu,\beta}$ such that
\[\|{\mathcal M}_f\|_{L^p(\T)}\leq C_{p,\nu,\beta}\|\tr\, f\|_{L^p(\T)}.\]
\item[(f)] Each $f\in H^p_\nu$ satisfies the maximum principle,
{\it i.e.} $|f|$ cannot assume a relative maximum in $\D$ unless 
it is constant. More generally, a non constant function in $H^p_\nu$
is open and discrete\footnote{A map is discrete if the preimage of
any value is a discrete subset of its domain.}.
\end{itemize}
\end{propositionsubsect}
\Rd It is rather easy to deduce Proposition \ref{cor:propHpnu} from the corresponding properties for the $G^p_{\alpha}(\D)$ class. Indeed, \Bk
by \eqref{kappa} and Proposition \ref{trick-hardy}
we can invert \eqref{correspw} by \eqref{correspf}, \Rd so that \Bk items $(a)$-$(e)$ follow at once from
Proposition \ref{gc-hardy}, Remark \ref{rmkY}, and the fact that
$f$ and $w$ share the same zeros because
$$
w=\frac{f-\nu\overline{f}}{\sqrt{1-\nu^2}}.
$$
To prove $(f)$, observe from Theorem \ref{trick2} that $w$, thus also
$f$ belongs to $W^{1,l}_{loc}(\D)$ for each $1<l<\infty$. Moreover
if we let \Rd $\nu_f(z):=\nu(z) \overline{\partial f(z)}/\partial f(z)$\Bk if 
$\partial f(z)\neq0$ and $\nu_f(z)=0$ otherwise, then $f$
is a pointwise a.e. solution in $\D$ of the \emph{classical} 
Beltrami equation:
\begin{equation}
\label{clBel}
\Rd \overline\partial f  = \nu_f  \, \partial f, ~~~~~~|\nu_f|\leq \kappa<1\Bk. 
\end{equation}
\Rd Indeed, $\left\vert \nu(z)\right\vert\leq \kappa$ for all $z\in \D$ and $\left\vert \overline{\partial f(z)}/\partial f(z)\right\vert=1$ when $\partial f(z)\neq 0$. \Bk It is then a standard result \cite[ Thm  11.1.2]{im2}
that $f=G(h(z))$, where $h$ is a quasi-conformal
topological map $\D\to\C$ satisfying \eqref{clBel} and $G$ a 
holomorphic function on $h(\D)$. The conclusion follows at once 
from the corresponding properties of holomorphic functions.
\hfill\fin\par

\begin{remarquesubsect}
\label{boundcnu}
When $\nu=0$, that is, when dealing with holomorphic Hardy spaces,
the best constant in \eqref{Fatouf} is $c_{0}=1$ 
because $\|f\|_{L^p(\T_r)}$ increases with $r$,
and then equality holds throughout. For general $\nu$,
a bound on $c_{\nu}$ depending solely on $\|\nu\|_{W^{1,\infty}(\D)}$
is easily derived from Proposition \ref{trick-hardy} and 
Theorem \ref{trick2}, but the authors do not know of a  sharp estimate.
\end{remarquesubsect}
\begin{remarquesubsect} \label{zero}
Assertion (c) in Proposition \ref{cor:propHpnu} implies that a 
function $f\in H^p_{\nu}(\D)$ whose trace is zero on a subset of $\T$ having  
positive Lebesgue measure must vanish identically.
\end{remarquesubsect}

As in the holomorphic case,
a function in $H^{p,0}_\nu$ is uniquely defined by its real part on $\T$:

\begin{propositionsubsect}
\label{prop:uniquenessf}
Let $f \in H^{p,0}_{\nu}$ be such that $\mbox{Re } (\tr\, f) = 0$ a.e. on
$\T$. Then $f \equiv 0$.
\end{propositionsubsect}
\Rd {\it Proof of Proposition \ref{prop:uniquenessf}.} 
Let $f\in H_\nu^{p,0}(\D)$ satisfy  $\mbox{Re}\,\tr\, f=0$ a.e. on $\T$. 
If we define $w$ through  \eqref{correspw}, then
$w\in G^{p,0}_\alpha$  by Proposition \ref{trick-hardy}
and clearly $\mbox{Re}\,\tr\, w=0$ a.e. on $\T$. Therefore 
$w\equiv0$ in view of Proposition \ref{gc-hardy}, assertion 4,
whence $f\equiv 0$ in $\D$. \hfill\fin\par \Bk

\noindent The next result shows that $H^p_\nu$ contains all $W^{1,p}(\D)$ 
solutions to (\ref{dbar}). That this inclusion is a strict one follows 
at once from Theorem \ref{thm:dirichletf} to come.
\begin{propositionsubsect} \label{hardysobol}
Let $f\in W^{1,p}(\D)$ be a solution to (\ref{dbar}). Then
$f\in H^p_{\nu}(\D)$, and there exists $C_{\nu,p}>0$ such that, 
\begin{equation}
\label{inegHSobt}
\left\Vert f\right\Vert_{H^p_{\nu}(\D)}\leq C_{\nu,p}\left\Vert f\right\Vert_{W^{1,p}(\D)}.
\end{equation}
\Rd Moreover, the trace of $f$ considered as an element of $W^{1,p}(\D)$ coincides with its trace seen as an element of $H^p_{\nu}(\D)$. 
\end{propositionsubsect}
Note that \eqref{inegHSobt} follows immediately from \eqref{Fatouf} and
the continuity of the trace operator from $W^{1,p}(\D)$ into $L^p(\T)$,
once it is known that $f\in H^p_\nu$.

The proof of Proposition \ref{hardysobol} {\Bl is}
given in Section \ref{proofpropf}.

\subsection{Regularity of the Dirichlet problem}
\label{lllllllll}
\subsubsection{\Rd Solvability of the Dirichlet problem in $G^p_{\alpha}(\D)$ \Bk}
\Rd We first focus on a slight variation of the Dirichlet problem for the $G^p_{\alpha}(\D)$ class. Let us introduce \Bk one more piece of notation by letting
\[P_+ \psi =\mbox{tr }\left({\mathcal C}\psi\right) 
\]
where the trace is a nontangential limit.
As is well-known, $P_+\psi$ exists a.e. on $\T$ as soon as $\psi\in L^1(\T)$,
but it may not lie in $L^1(\T)$. If, however, $1<p<\infty$,
then $P_+$ is a continuous projection from $L^p(\T)$ onto $\tr\,H^p$
called the analytic projection \cite[Ch. III, Sec. 1]{gar}.
It is an interesting variant of the Dirichlet problem to
solve equation (\ref{eq:w}) while prescribing the analytic projection of 
the solution on $\T$. As the next theorem shows, 
$G^{p}_{\alpha}$ is a natural space for this.
\begin{theoremsubsubsect}
\label{dirichletp+}
For $\alpha\in L^\infty(\D)$ and $g\in {H}^{p}$, 
there is a unique $w\in G^p_{\alpha}$ such that
\begin{equation} \label{p+w}
P_{+}(\mbox{tr }w)= \tr\, g\,.
\end{equation}
This solution satisfies
\begin{equation} \label{formula}
w=g+T(\alpha \overline{w}) \, ,\ \mbox{a.e. in } \D\,,
\end{equation}
 and it holds that
\begin{equation} \label{estimwphi}
\left\Vert w\right\Vert_{G^p_{\alpha}}\leq C_{p,\alpha} 
\left\Vert g\right\Vert_{H^p(\D)}. 
\end{equation}
\end{theoremsubsubsect}
\Rd Here is now the solution of the (usual) Dirichlet problem for the class
$G^p_{\alpha}$: \Bk
\begin{theoremsubsubsect}
\label{thm:dirichletw}
Let $\alpha\in L^\infty(\D)$ and $\psi \in L^p_\RR(\T)$. 
\begin{itemize}
\item[(a)] To each $c\in\R$,  there uniquely exists 
$w \in G^{p}_{\alpha}$  such that
$\mbox {Re } (\tr\, w) = \psi$  a.e. on $ \T$ and 
$\int_0^{2\pi}\mbox{Im}\,\tr\,
w(e^{i\theta})\,d\theta=2\pi c$. Moreover there are constants $c_{p,\alpha}$
and $c_{p,\alpha}'$ such that
\begin{equation}
\label{inegDw}
\left\Vert \tr\,  w\right\Vert_{L^p(\TT)} \leq
c_{p, \alpha} \left\Vert \psi \right\Vert_{L^p(\TT)}+
c_{p, \alpha}'\, |c|\,.
\end{equation}
\item[(b)] When \eqref{alphabeta} holds, there uniquely exists 
$w \in G^{p,0}_{\alpha}$ such that
$\mbox {Re } (\tr\, w) = \psi$  a.e. on $ \T$. Furthermore, there is a constant
$c_{p,\alpha}''$ such that
\begin{equation}
\label{inegDw0}
\left\Vert \tr\,  w\right\Vert_{L^p(\TT)} \leq
c_{p, \alpha}'' \left\Vert \psi \right\Vert_{L^p(\TT)} \, .
\end{equation}
\end{itemize}
\end{theoremsubsubsect}

The proofs of Theorems \ref{dirichletp+}-\ref{thm:dirichletw} are
given in section \ref{proofdirw}.

\subsubsection{\Rd Solvability of the Dirichlet problem in $H^p_{\nu}(\D)$ \Bk}
The following \Rd result \Bk shows that
$H^p_\nu$ is the natural space to consider when handling $L^p$ 
boundary data in \eqref{dbar} and \eqref{div}.
\begin{theoremsubsubsect}
\label{thm:dirichletf}
For all $\varphi \in L^p_\RR(\T)$, there uniquely exists 
$f \in H^{p,0}_{\nu}$ such that, a.e. on $\T$:
\begin{equation}
\mbox{Re } (\tr\, f) = \varphi \, .
\label{Rephi}
\end{equation}
Moreover, there exists $c_{p, \nu} > 0$ such that:
\begin{equation}
\label{inegRef}
\left\Vert f\right\Vert_{H^p_{\nu}(\D)}\ \leq
c_{p, \nu} \left\Vert \varphi \right\Vert_{L^p(\TT)} \, .
\end{equation}
\end{theoremsubsubsect}
\Rd From Proposition \ref{hardysobol}, Theorem \ref{thm:dirichletf} clearly extends Theorem \ref{dirichlet} when the boundary data belong to $L^p(\T)$. \Bk\par

\medskip

\noindent \Rd Let us give at once the proof of Theorem \ref{thm:dirichletf}, which is quite easy to deduce from previous statements. \Bk 
Define $\alpha$ through \eqref{alphabeta} and put
$\psi = \varphi \sigma^{1/2} \in L^p_\RR(\T)$. Apply Theorem
\ref{thm:dirichletw}, point (b), to obtain  $w
\in G^{p,0}_{\alpha}$ such that $\mbox{Re } (\tr\, w) = \psi$. 
If we let $f$ be given  by (\ref{correspf}),
then  $f\in H^{p,0}_\nu$ by 
Proposition \ref{trick-hardy}, point (a). Moreover, from \eqref{correspw},
we see that \eqref{Rephi} holds. The uniqueness of $f$ comes from
Proposition \ref{prop:uniquenessf}.
Inequality \eqref{inegRef}  follows from \eqref{correspf},
\eqref{ellipticsigma}, \eqref{inegDw0} and \eqref{Fatou}.
\hfill\fin\par

\medskip

\noindent Dwelling on Proposition
\ref{cor:propHpnu} and Theorem \ref{thm:dirichletf}, we are now able to 
derive an analog of the Fatou theory 
\cite[Ch. I, Sec. 5]{gar} for \eqref{div}, at least when
$1<p<\infty$. It should be compared to classical results on the Dirichlet 
problem in Sobolev classes \cite{campanato,gt}. For once, we recall all the 
assumptions to ease this comparison.
\begin{theoremsubsubsect}
\label{thm:dirichletu} Let $1<p<\infty$ and $\sigma\in W^{1,\infty}_\R(\D)$ 
satisfy \eqref{ellipticsigma}. Any $u\in L^p_\R(\D)$ satisfying \eqref{div} in the sense of distributions 
lies in $W^{1,l}_{\R,loc}(\D)$ for all $l\in(1,\infty)$. If, moreover,
\begin{equation}
\label{Hardypu}
\|u\|_{F,p}:=\supess_{0<r<1}
 \left\Vert u\right\Vert_{L^p(\T_r)} <+\infty \, ,
\end{equation}
then $u$ has a nontangential limit $\tr\, u$ a.e. on $\T$ which is also the 
limit of $e^{i\theta}\mapsto u(re^{i\theta})$ in $L^p_\R(\T)$
as $r\to 1^-$. In this case, for ${\cal M}_{|u|}$
the non tangential maximal function, we have
\[\|\tr\, u\|_{L^p_\RR(\T)}\leq \|u\|_{F,p}\leq 
c_{p,\nu}\|{\cal M}_{|u|}\|_{L^p(\T)}\leq 
C_{p,\nu} \|\tr\, u\|_{L^p_\RR(\T)}.\]
Conversely, each member of $L^p_\R(\T)$ is uniquely the non tangential limit 
of some distributional solution $u\in L^p(\D)$ of \eqref{div} 
satisfying $\|u\|_{F,p}<+\infty$.
\end{theoremsubsubsect}
\Rd The proof  of Theorem  \ref{thm:dirichletu}
is carried out in  Section \ref{proofdirf}. \Bk

Theorem \ref{thm:dirichletf} allows one to define a
generalized conjugation operator $\cH_\nu$ from $L^p(\TT)$ into itself,
that was introduced on $W^{1/2,2}(\T)$ in \cite{ap2}  {\Bl as the 
$\nu$-Hilbert transform}.
More precisely, to each $\varphi\in L^p_\R(\T)$, we associate  
the unique function $f\in H^{p,0}_\nu$ such that
Re tr $f=\varphi$,  and we set 
${\cH}_\nu\varphi= \mbox{Im} \,tr f\in L^p(\T)$. 
It now follows from Theorems 
\ref{dirichlet} and \ref{thm:dirichletf} that:
\begin{corollarysubsubsect}
\label{MRieszgen}
The operator ${\cH}_{\nu}$ is bounded both on
$L^p_{\R}(\T)$ and on $W^{1-1/p,p}_{\R
}(\T)$.
\end{corollarysubsubsect} 
When $\nu=0$, we observe that $ {\cH}_0 \, \varphi$ is just the
harmonic conjugate\footnote{\Rd 
{\Bl Though it has the same behaviour, it is distinct from the Hilbert transform, see \cite[Chap. III, So. 1]{gar}}.\Bk} of $\varphi$ normalized to have zero mean on $\T$.
That the operator  ${\cH}_0$ is continuous from 
$L^p_\R(\T)$ into itself is the well-known
M. Riesz theorem \cite[Ch. III, thm. 2.3]{gar}. 
Corollary \ref{MRieszgen} thus generalizes the latter.

\subsubsection{\Rd Improved regularity results for the Dirichlet problem \Bk}
We turn to higher regularity for solutions to (\ref{dbar}). More precisely, 
we shall study the improvement in the conclusion of Theorem
\ref{thm:dirichletf} when the boundary condition lies in
$W_\R^{1,p}(\T)\subset W_\R^{1-1/p,p}(\T)$.
First, the generalized conjugation operator preserves this smoothness class
(compare Corollary \ref{MRieszgen}): 
\begin{propositionsubsubsect} \label{hnusob}
The operator ${\mathcal H}_{\nu}$ is bounded on $W^{1,p}_{\R}(\T)$.
\end{propositionsubsubsect}
Next, assuming in Theorem \ref{thm:dirichletf} that 
$\varphi\in W^{1,p}_{\R}(\T)$, not only does $f$ belong to $W^{1,p}(\D)$,
as predicted by Theorem \ref{dirichlet}, but 
the derivatives of $f$ satisfy a condition of Hardy type: 
\begin{theoremsubsubsect} \label{sobnontang}
Let $\varphi\in W_{\R}^{1,p}(\T)$ and $f\in W^{1,p}(\D)$ be the
unique solution to \eqref{dbar} on $\D$  satisfying
${\rm Re} \, (\tr\, f)=\varphi$ and such that (\ref{normal:f}) holds. Then,
\begin{itemize}
\item[$(a)$]
$\tr\, f\in W^{1,p}(\T)$, and it holds that
\begin{equation}
\label{boundineg}
\|\tr\, f\|_{W^{1,p}(\T)}\leq C_{ p, \nu} \|\varphi\|_{W^{1,p}(\T)}.
\end{equation}
\item[$(b)$]
The functions $\partial f$ and $\overline{\partial} f$
satisfy a Hardy condition of the form
\begin{equation}
\label{HCp1}
 \supess_{0<r<1}
 \left\Vert \partial f\right\Vert_{L^p(\T_r)}\leq C_{ p, \nu}
\|\tr\, f\|_{W^{1,p}(\T)},
\end{equation}
\begin{equation}
\label{HCp2}
 \supess_{0<r<1}
 \left\Vert \overline{\partial} f\right\Vert_{L^p(\T_r)}\leq C_{ p, \nu}
\|\tr\, f\|_{W^{1,p}(\T)},
\end{equation}
and for the non tangential maximal function of $\|\nabla f\|$,
it holds that
\begin{equation}
\label{inegmaxgrad}
\|{\mathcal M}_{\|\nabla f\|}\|_{L^p(\T)}
\leq C_{p,\nu,\beta}\|\tr\, f\|_{W^{1,p}(\T)},
\end{equation}
where $\nabla f(\xi)\in\C^2$ 
is the gradient of $f$ and $\beta$ the aperture of 
the sectors in (\ref{defntmax}).
\item[$(c)$]
If we define $\Phi \in L^p(\T)$ by
\begin{equation}
\label{relim}
\Phi (e^{i\theta}):=-ie^{-i\theta}\,\,\frac{\partial_{\theta}(\tr\, f)(e^{i\theta})-
\nu(e^{i\theta})\partial_{\theta}\overline{(\tr\, f)(e^{i\theta})}}{1-\nu^2(e^{i\theta})},
\end{equation}
then $\partial f$ and $\overline{\partial} f$ have non tangential 
limit $\Phi$ and $\nu\overline{\Phi}$
a.e. on $\T$, and $\partial f(r e^{i\theta})$,
$\overline{\partial} f(r e^{i\theta})$ converge in $L^p(\T)$ to their
respective nontangential limits as $r\to1$.
\end{itemize}
\label{thmLB}
\end{theoremsubsubsect}
Clearly, (a) is a rephrasing of
Proposition \ref{hnusob}. \par

\medskip

\noindent As a corollary of Theorem \ref{sobnontang}, we obtain the following result 
(compare \cite{FJR}), which plays for the Neumann problem the same role as Theorem \ref{thm:dirichletu} does for the Dirichlet problem:
\begin{corollarysubsubsect} \label{fjr}
\begin{itemize}
\item[$1.$]
Let $u\in W^{1,p}_\R(\D)$ be a solution of $\mbox{div }(\sigma\nabla u)=0$ in $\D$ such that $\nabla u$ satisfies the following Hardy condition:
\begin{equation} \label{formhardy}
\supess_{0<r<1} \left\Vert  \nabla u \right\Vert_{L^p(\T_r)}<+\infty.
\end{equation}
Then $\tr u\in W^{1,p}(\T)$ and $u\in W^{1,p_1}_\R(\D)$ for every $p_1 \in (p, 2p)$. Moreover, ${\mathcal M}_{\|\nabla u\|}\in L^p(\T)$, and there exists a vector field $\Phi\in L^p(\T,\R^n)$ such that $\nabla u\rightarrow \Phi$ n.t. almost everywhere on $\T$. In particular, $\partial_nu\in L^p(\T)$, and one has $\int_{\T} \sigma\partial_nu=0$.
\item[$2.$]
Conversely, if $g\in L^p_{\R}(\T)$ satisfies $\int_{\T} \sigma g=0$, there exists a function $u\in W^{1,p}_\R(\D)$ solving $\mbox{div }(\sigma\nabla u)=0$ in $\D$ such that $\nabla u$ satisfies a Hardy condition of the form (\ref{formhardy}),  ${\mathcal M}_{\|\nabla u\|}\in L^p(\T)$ and $\partial_nu=g$ on $\T$. Moreover, $u$ is unique up to an additive constant. 
\end{itemize}
\end{corollarysubsubsect}
\Bk
All these results will be established in 
Section \ref{proofhigherreg}. 
\subsection{Density of traces}
\label{mmmmmmmmmmmmmm}
We come to some density properties of
traces of solutions to (\ref{dbar}). Loosely speaking, they assert that 
if  $E\subset \T$ is not too large, then every complex function on $E$ 
can be approximated by the trace of a solution to \eqref{dbar} on $\D$.

\subsubsection{\Rd Density in Sobolev spaces \Bk}
\noindent We say that an open subset $I$ of $T$ has the 
\emph{extension property} if every function in $W^{1,p}(I)$ is 
the restriction to $I$ of some function in $W^{1,p}(\T)$. If $I$ is a proper
open subset of $\T$, it decomposes into a countable union of disjoint open 
arcs $(a_j,b_j)$ and the extension property is equivalent to the fact that
no $a_j$ (resp. $b_j$) is a limit point of the sequence $(b_k)$   (resp. $(a_k)$).

We begin with a density property of Sobolev solutions to \eqref{dbar}
on proper extension subsets:
\begin{theoremsubsubsect}\label{densiteW}
Let $I\neq \T$ be an open subset of $T$ having the extension property.
Then, the restrictions to $I$ of traces of
$W^{1,p}(\D)$-solutions to \eqref{dbar}
form a dense subspace of $W^{1-1/p,p}(I)$.
\end{theoremsubsubsect}
This should be held in contrast with the fact that
the traces on $\T$ of $W^{1,p}(\D)$-solutions to (\ref{dbar}) 
form a proper closed subspace of $W^{1-1/p,p}(\T)$. 

The proof of Theorem \ref{densiteW} is given in Section \ref{proofdensesob}.

\subsubsection{\Rd Density in Lebesgue spaces\Bk}
By the density of $W^{1-1/p,p}(I)$ in $L^p(I)$, Theorem \ref{densiteW}
easily implies that $(\tr H^p_\nu)_{|_I}$ is a dense subset of 
$L^p(I)$ for $I$ a proper 
open subset of $\T$ having the extension property. The fact that this remains 
true as soon as $I$ is not of full measure lies a little deeper:
\begin{theoremsubsubsect}
\label{thm:densite}
Let $I\subset\T$ be a measurable subset such that $\T \setminus I$ has positive Lebesgue measure.
The restrictions to $I$
of traces of $H_\nu^p$-functions are dense in $L^{p}(I)$. 
\end{theoremsubsubsect}
\begin{remarquesubsubsect}
\label{rmkinf}
When $I\subset\T$ is not of full measure and $\phi\in L^p(I)$, 
Theorem \ref{thm:densite} entails there is a sequence of functions 
$f_k\in H^p_\nu$ whose trace on $I$ converges to $\phi$ in
$L^{p}(I)$. Now, since  balls in $\tr\, H^p_\nu$ are weakly compact
by Proposition \ref{cor:propHpnu} point (b),
it must be that either $\phi$ is the trace on $I$ of a 
$H^p_\nu$-function or 
$\left\Vert \tr\, f_k\right\Vert_{L^{p}(\T\setminus I)}\rightarrow
+\infty$ with $k$. In view of Theorem \ref{densiteW}, the corresponding remark 
applies when $I$ is an open subset of $\T$ with the extension property 
and $\varphi\in W^{1-1/p,p}(I)$ gets approximated in this space
by a sequence of traces of $W^{1,p}(\D)$-solutions to \eqref{dbar}.
\end{remarquesubsubsect}
It is worth recasting Remark \ref{rmkinf} in terms of ill-posedness of the
inverse Dirichlet-Neumann problem from incomplete boundary data. Indeed,
assume that $u$ satisfies \eqref{div} and, say, $\tr\, u\in W^{1,p}(\T)$.
Observe that the normal derivative $\partial_n u$ exists as a 
nontangential limit 
in $L^p(\T)$ by Theorem \ref{thmLB}. Thus, upon rewriting 
\eqref{system} on $\T$ in the form 
\begin{equation}
\label{CRtang}
\left\{
\begin{array}{l}
\partial_nv=-\sigma\partial_\theta u,\\
\partial_\theta v=\sigma\partial_n u,
\end{array}
\right.
\end{equation}
we see that the knowledge of $\tr\,u$ and $\tr\,\partial_n u$ 
on some arc
$I\subset \T$
is equivalent  to the knowledge 
on $I$ of $\tr\, f$ where $f=u+iv$ meets \eqref{dbar} 
with boundary conditions $\mbox{Re} f=u$ and, say, $\int_I v=0$. 
Note from Proposition \ref{cor:propHpnu} that this determines $f$ completely.
Now, if the knowledge of $\tr\,f$ gets corrupted by 
measurements and rounding off errors, as is the case in computational and 
engineering practice, it can still be approximated arbitrarily well in
$W^{1-1/p,p}(I)$ by a solution to \eqref{dbar} but the trace of the 
latter will 
grow large in $L^p(\T\setminus I)$, {\it a fortiori} in 
$W^{1-1/p,p}(\T\setminus I)$ when the approximation error gets 
small.

The proof of Theorem \ref{thm:densite}
 is carried out in Section \ref{proofdensehard}.
\subsection{Duality} \label{dual}
Keeping in mind that $1<p<\infty$ and $1/p+1/q=1$, 
we introduce a duality pairing on 
$L^p(\T)\times L^q(\T)$, viewed as real vector spaces,
 by the formula:
\begin{equation}
\label{dualpair}
\langle f , g \rangle = \mbox{Re} \, \int_0^{2 \pi} f \, g \, 
\frac{d \theta}{2\pi}.
\end{equation}
Clearly this pairing isometrically identifies $L^q(\TT)$ with 
the dual of $L^p(\TT)$. The fact that $H^p$ is the orthogonal space to
$e^{i\theta} H^q$ under \eqref{dualpair} is basic to the dual 
approach of extremal problems
in holomorphic Hardy spaces \cite[Ch. 8]{duren}. In this section,
we derive the corresponding results for the spaces $H^p_\nu$. Recall that
$\partial_t=\partial_\theta/2\pi$ on $\T$.
 \begin{propositionsubsect}
 \label{orthogformula}
The orthogonal to $\tr\, H^p_\nu$ under the duality pairing defined in \eqref
{dualpair} is
\[
(\tr\, \, H^{p}_\nu)^\perp=\partial_\theta \, \left(\tr\, H^{q}_{-\nu} \cap W^{1,q}(\T) \right) \ . 
\]
\end{propositionsubsect}
Proposition \ref{orthogformula} and the Hahn-Banach theorem now
team up to  yield:
\begin{theoremsubsect}
\label{thm:dualite}
\mbox{} 
\begin{itemize}
\item[(i)] Under the pairing \eqref{dualpair}, the dual space $(\tr\, H^p_\nu)^*$ 
of $\tr\, H^p_\nu$ is naturally
isometric to the quotient space $L^q(\T)/(\tr\, H^p_\nu)^\perp$, 
that is
\[(\tr\, H^p_\nu)^* \sim
L^q(\T)/ \left(\partial_\theta \left( \tr\, H^{q}_{-\nu} \cap W^{1,q}(\T) \right)\right) \, .
\]
\item[(ii)] For each $\Phi\in L^q(\T)$, it holds the duality relation
\begin{equation}
\label{dualnue}
\inf_{g\in \partial_\theta (\tr\, H^q_{-\nu}\cap W^{1,q}(\T))}
\|\Phi-g\|_{L^q(\T)}=\sup_{\stackrel{f\in H^p_\nu}{\|\tr\, f\|_{L^p(\T)=1}}}
\frac{1}{2\pi}\mbox{{\rm Re}} \, \int_0^{2 \pi} \Phi \, \tr\, f \, d \theta.
\end{equation}
\item[(iii)] For each $\Psi\in L^p(\TT)$, it holds the duality relation
\begin{equation}
\label{dualnue1}
\inf_{f\in H^p_\nu}\|\Psi-\tr\, f\|_{L^p(\T)}=
\sup_{\stackrel{g\in \tr\, H^q_{-\nu}\cap W^{1,q}(\T)}
{\|\partial_\theta g\|_{L^q(\T)}=1}}
\frac{1}{2\pi}\mbox{{\rm Re}} \, \int_0^{2 \pi} \Psi \, \partial_\theta g \, d \theta.
\end{equation}
\end{itemize}
\end{theoremsubsect}
\Rd Granted Proposition
\ref{orthogformula}, Theorem \ref{thm:dualite} is a standard application of the Hahn-Banach 
theorem \cite[{ Thms} 7.1, 7.2]{duren}. \Bk\hfill\fin
\begin{remarquesubsect}
It is easy to check (compare \cite[{ Thm} 3.11]{duren}) that 
$\partial_\theta(\tr\, H^q\cap W^{1,q}(\T))=e^{i\theta}H^q$, hence 
\eqref{dualnue}-\eqref{dualnue1} reduce to standard duality relations
in  Hardy spaces when $\nu=0$.
\end{remarquesubsect}
The proofs of Proposition \ref{orthogformula} and Theorem \ref{thm:dualite} will be given in Section \ref{sec:charac}.
\section{Proofs} \label{proofmain}

\subsection{\Rd The Dirichlet problem in Sobolev spaces \Bk}
Let us give first the proof of Theorem \ref{dirichlet}. Put $\sigma=(1-\nu)/(1+\nu)$ , so that
$\sigma\in W^{1,\infty}(\D)$ satisfies \eqref{ellipticsigma}. 
By \cite{campanato}, there uniquely exists 
$u\in W^{1,p}(\D)$ meeting 
$\tr\, u=\varphi$ for which \eqref{weaksol} holds with $\Omega=\D$;
moreover, by the open mapping theorem, one has
\begin{equation}
\label{campu}
\|u\|_{W^{1,p}(\D)}\leq c_{p,\nu} \left\Vert \varphi\right
\Vert_{W^{1-1/p,p}(\T)}.
\end{equation}
Put
\begin{equation}
\label{defHodge}
{\mathcal G}_{0,p}:=\{\nabla g;\ g\in W^{1,p}_{0,\R}(\D)\}
\ \ \ \mbox{and}\ \ \ 
{\mathcal D}_q:=\{(\partial_y h,-\partial_x h)^T;\ h\in W^{1,q}_\R(\D)\}.
\end{equation}
Proceeding by density on the divergence formula for smooth functions, 
we easily get
\begin{equation}
\label{orthogHodge}
\langle G,D\rangle:=\int_\D G.D\,dm=0,\ \ \ \ G\in{\mathcal G}_{0,p},\ D\in{\mathcal D}_q.
\end{equation}
Now, by Hodge theory \cite[{ Thm} 10.5.1]{im2}\footnote{The result is stated
there using the language of differential forms that
we did not introduce here.}, each vector field in $L^p(\D)\times L^p(\D)$ 
(resp. $L^q(\D)\times L^q(\D)$) is uniquely the
sum of a member of ${\mathcal G}_{0,p}$ (resp. ${\mathcal G}_{0,q}$) and a 
member of  ${\mathcal D}_p$ (resp. ${\mathcal D}_q$).
If we set accordingly $\sigma\nabla u=G+D$, we gather by density from 
\eqref{weaksol} and \eqref{orthogHodge}
that $\langle G, V\rangle=0$ for every $V\in L^q(\D)\times L^q(\D)$, 
implying that
$\sigma\nabla u\in{\mathcal D}_p$. In other words there is 
$v\in W^{1,p}_{\R}(\D)$ for which \eqref{system} holds, thus
by inspection $f=u+iv\in W^{1,p}(\D)$ satisfies
(\ref{dbar}) pointwise a.e. on $\D$. Then, $f$ satisfies \eqref{dbar}
in the distributional sense as well.
To check $f$ is unique, subject to
$\mbox{Re}\,\tr\, f=\varphi$,  up to
an additive pure imaginary constant, observe if $f\in W^{1,p}(\D)$ satifies
\eqref{dbar} that $u=\mbox{Re}f$ and $v=\mbox{Im} f$ both lie in
$W^{1,p}_\R(\D)$ and that \eqref{system} holds. Therefore 
$\sigma\nabla u \in {\mathcal D}_p$ and, in view of \eqref{orthogHodge},
we see that \eqref{weaksol} holds with $\Omega=\D$. As such a $u$ is uniquely
defined by $\tr\, u=\varphi$, we conclude that
 $v$ is uniquely defined by \eqref{system}, up to an additive constant.

Next,we  observe from (\ref{system}) and (\ref{ellipticsigma})
that $\|\nabla v\|_{L^p(\D)}\leq C \|\nabla u\|_{L^p(\D)}$.
Therefore, if $v$ gets normalized by \eqref{normal:f}, it follows from 
\eqref{campu} and the Poincar\'e inequality 
\cite[Ch. 4, Ex. 4.10]{Ziemer} that 
\[\|v\|_{W^{1,p}(\D)}\leq C_p \|\nabla v\|_{L^p(\D)}\leq C
C_pc_{p,\nu} \left\Vert \varphi\right
\Vert_{W^{1-1/p,p}(\T)}\] so that (\ref{estimf}) indeed holds. 
 \hfill\fin\par

\subsection{Preliminaries on spaces and operators} \label{op}
In the present subsection, we recall some properties of the operators 
${\mathcal C}$ and $T$, introduced in Section \ref{gpalphaprop}, and of 
the Beurling operator appearing in equation \eqref{Beurldef} below.
\par

For $h\in L^p(\C)$, we define the operator
$\breve{T}$ by
\[
\breve{T}h(z)=\frac1{2\pi i}\iint_{\D}
\frac{h(\xi)}{\xi-z}d\xi\wedge
d\overline{\xi}\,,\  \ z\in \C\,.
\]
If $w\in L^p(\D)$ and $\breve{w}$ is the extension of $w$ by 0 outside
$\D$, then obviously $(\breve{T} \breve{w})_{|_\D} = T w$.
Next, for $u\in L^p(\C)$, we denote by $S$ the Beurling operator:
\begin{equation}
\label{Beurldef}
Su(z)=\lim_{\e\to 0+}\frac1{2\pi i} \iint_{\xi\in \C,\ |\xi-z|\geq\e}
\frac{u(\xi)}{(\xi-z)^{2}}d\xi\wedge d\overline{\xi}
\,,\  \ \mbox{a.e.}\ z\in \C\,.
\end{equation}
The existence of $Su$ a.e. follows from the Calder\`on-Zygmund
theory of singular integral operators 
\cite[Ch. II, { Thm} 4]{stein}.
Here are the
properties of ${\mathcal C}$, $\breve T$, $T$ and $S$ that we use:
\begin{propositionsubsect} \label{properties}
Let as usual $1<p<+\infty$. Then the following assertions hold.
\begin{itemize}
\item[$1.$] The Cauchy operator ${\mathcal C}$ is bounded from
$L^p(\T)$ onto $H^p(\D)$ and from $W^{1-1/p,p}(\T)$ to $W^{1,p}(\D)$.
\item[$2.$] The Beurling operator $S$ is bounded from $L^p(\C)$ into
itself. 
\item[$3.$] The operator $\breve{T}$ maps 
$L^p(\C)$ continuously into $W^{1,p}_{loc}(\C)$.
\item[$4.$] The operator
$T$ is bounded from $L^p(\D)$ into $W^{1,p}(\D)$, 
and is compact from $L^p(\D)$ to $L^p(\D)$. Moreover
$\bar\partial T w =w$ and $\partial T w =(S\breve{w})_{|_\D}$
for all $w\in L^p(\D)$.
For any  $\alpha\in L^{\infty}(\D)$ the operator 
$w\mapsto w-T(\alpha\overline{w})$ is an isomorphism of $L^p(\D)$.
\end{itemize}
\end{propositionsubsect}
The next result will be of technical importance to establish the regularity 
properties of $G_\alpha^p$-functions, compare Remark \ref{rmkY}.
\begin{lemmasubsect} \label{improve}
Let $p\in (1,+\infty)$ as always, and $\alpha\in L^\infty(\D)$.
\begin{itemize}
\item[$1.$]
If $g\in H^p(\D)$, then $g\in L^{p_1}(\D)$ for 
$p_1\in [p,2p)$.
\item[$2.$]
If $w \in L^p(\D)$ and if $w - T(\alpha\overline{w}) \in H^p$, then 
there is $p^* > 2$ such that
$T(\alpha\overline{w}) \in W^{1,{p^*}}(\D) \subset C^{0,1-2/p^*}(\overline{\D})$ and $\breve{T}\left(\breve{\alpha \overline{w}}\right) \in W^{1,{p^*}}_{loc}(\C) \subset C_{loc}^{0,1-2/p^*}(\C)$. Moreover, 
\begin{equation}\label{td}
\|T(\alpha \overline{w})\|_{W^{1,p^*}(\D)}\leq C_{p,\alpha}\left(\|w\|_{L^p(\D)}+\|w-T(\alpha \overline{w})\|_{H^p(\D)}\right).
\end{equation}
\end{itemize}
\end{lemmasubsect}
\Rd In order not to disrupt the reading, we postpone the proofs of Proposition \ref{properties} and Lemma \ref{improve} to Appendix \ref{technic}. \Bk
\subsection{Factorization and boundary behaviour  in $G^p_{\alpha}$}\label{secprop:w}

This section is devoted to the proof of Theorem \ref{trick2} and 
Proposition \ref{gc-hardy}. In the proof of the former, we make use of the following Lemma.
\begin{lemmasubsect} \label{w1qloc}
Let $r\in L^{\infty}(\C)$ be supported in $\D$. Then, the function
$$
u(z)=\iint_{\D} \frac{z\overline{r(\zeta)}}
{1-\bar\zeta z} d\zeta\wedge d\bar\zeta\, ,\ \ \ z\in\C\,,
$$
is holomorphic in $\D$ and 
belongs to $W^{1,l}_{loc}(\C)$ for all $l\in (1,+\infty)$.
\end{lemmasubsect}
{\it Proof of Lemma \ref{w1qloc}.}  The function  $u$ is clearly holomorphic in $\D$,
{\it a fortiori} $u\in W^{1,l}_{loc}(\D)$. 
It is therefore enough to show, say, that 
$u\in W^{1,l}_{loc}(\C\setminus  \overline{\D}_{1/2} )$
for all $l\in(1,+\infty)$. In turn, it is sufficient to prove that
$$
u(1/\overline{z})=-
\overline{\int_{\D} \frac{r(\zeta)}{z-\zeta}d\zeta\wedge d\bar\zeta},
$$
lies in $W^{1,l}(\D_2)$. The conclusion now follows from
assertion $3$ in Proposition \ref{properties}.
\hfill\fin\par

\medskip

{\it Proof of Theorem \ref{trick2}.} put
$r(z)=\alpha(z)\overline{w(z)}/w(z)$ if $w(z)\not=0$ and $r(z)=0$ if 
$w(z)=0$ or $z\not\in\D$. Then $r\in L^\infty(\C)$ and 
$\|r\|_{L^\infty(\D)} \leq\|\alpha\|_{L^\infty(\D)}$.
Define 
\begin{equation}
\label{defs}
s(z)=\frac1{2\pi i}\iint_\D
\left(\frac{r(\zeta)}{\zeta-z}+\frac{z\overline{r(\zeta)}}
{1-\bar\zeta z}\right)
d\zeta\wedge d\bar\zeta \, , \mbox{ for } z \in \C \, .
\end{equation}
and observe, from Lemma \ref{w1qloc}
and assertion $3$ in Proposition \ref{properties}, that 
$s\in W^{1,l}_{loc}(\C)$ for all $l\in(1,+\infty)$.
In particular $s$ is continuous and,  since
$1/z=\overline{z}$ for $z\in \T$, we see from \eqref{defs}
that $\mbox{Im}\,s(z)=0$ there.
Also, assertion 4 of Proposition \ref{properties} and Lemma \ref{w1qloc}
show that $\bar\d  s=r$ in $\D$. Furthermore, a straightforward majorization
gives us for $z\in\C$ 
\begin{equation}
\label{majos}
|s(z)|\leq\frac{\|\alpha\|_{L^\infty(\D)}}\pi\iint_\D
\left(\frac1{|\zeta-z|}+\frac1{|\zeta-{1/z}|}\right) 
\,dm \leq 4\|\alpha\|_{L^\infty(\D)} \, ,
\end{equation}
thus \eqref{estims} holds.
Next, we put $F=e^{-s}w$ and claim that $F$ is holomorphic in $\D$.
Indeed, $F\in L^p(\D)$ hence,
by Weyl's lemma \cite[Thm 24.9]{Forster}, it is enough to check 
that $\bar\d F=0$ on $\D$ in the sense of distributions. 
Let $\psi\in \cD(\D)$ and $\psi_n$ a sequence in $\cD(\R^2)_{|_\D}$ 
converging to $s$ in $W^{1,l}(\D)$ for some $l>\max(q,2)$. Thus $\psi_n$ 
converges boundedly to $s$ in $W^{1,q}(\D)$ by the Sobolev imbedding theorem 
and H\"older's inequality. Then, by dominated convergence,
$$
\langle \bar\d F,\psi\rangle=
-\langle e^{-s}w,\bar\d\psi\rangle=
-\lim_n\langle w, e^{-\psi_n} \bar\d\psi\rangle=
-\lim_n\langle w,\bar\d(e^{-\psi_n} \psi)+\psi e^{-\psi_n}\bar\d\psi_n\rangle$$
$$
=\lim_n\langle \alpha\overline{w}, e^{-\psi_n}\psi\rangle
-\lim_n\langle w,\psi e^{-\psi_n}\bar\d\psi_n\rangle
=\langle e^{-s}(\alpha\overline{w}-w\bar\d s),\psi\rangle=0
$$
since $w\bar\d s=wr=\alpha \overline w$, where we used in the fourth 
equality that $e^{-\psi_n}\psi\in\cD(\D)$. This proves the claim
and provides us with \eqref{factBN} where \Rd $\mbox{Im}\,\tr\, s=0$. \Bk
Now, by the Sobolev imbedding theorem, $s$ is bounded, and since
$\exp$ is locally Lipschitz on $\C$
it follows that $w\in W^{1,l}_{loc}(\D)$ for all  $l\in(1,\infty)$.
Finally, by the boundedness of $s$, it is clear from the definitions
that $w\in G_\alpha^p$ if, and only it $F\in H^p$.
In this case, it 
follows from Lemma \ref{improve}, point 1, that $w = e^s \, F \in L^{p_1}(\D)$ for all $p_1 \in [p, 2p)$.

To obtain from \eqref{factBN} another factorization $w=e^{s_1}F_1$, 
where this time
$\mbox{Re}\,\tr\, s_1=0$, it is enough to change the ``$+$'' sign into a ``$-$'' one in the definition \eqref{defs} of $s$ .
\hfill\fin\par

\medskip

\noindent For the proof of Proposition \ref{gc-hardy}, we need the following 
version of the Cauchy-Green formula.
\begin{lemmasubsect} \label{gc}
When $\psi\in W^{1,p}(\D)$, it holds for almost every $z\in \D$ that
\begin{equation} \label{green}
\psi(z)={\mathcal C}\left( \tr \, \psi \right)(z)+\frac1{2\pi i}
\iint_{\D} 
\frac{\overline{\partial}\psi(\xi)}{\xi-z}d \, \xi \wedge d \, \overline{\xi}.
\end{equation}
\end{lemmasubsect}
{\it Proof.} 
Note that (\ref{green}) means
$\psi={\cal C}(\tr\, \psi)+T(\bar\d\psi)$. For $\psi\in\cD(\RR^2)$, 
this is standard \cite[thm. 1.2.1]{Hormander}. In general $\psi$ is the 
limit in $W^{1,p}(\D)$ of a sequence 
$(\psi_{n})_{n\in \NN}\in \left.{\mathcal D}(\RR^{2})\right\vert_{\D}$.
By continuity of the trace and Proposition \ref{properties},
items 1, 4, the conclusion follows from taking a pointwise convergent 
subsequence of the $L^p(\D)$ convergent sequence $T(\bar\d\psi_n)$. 
\hfill\fin\par
\medskip
{\it Proof of Proposition \ref{gc-hardy}.} Let $w\in G_\alpha^p$. 
By Theorem \ref{trick2}, we have $w=e^s F$ where $s\in W^{1,l}(\D)$,
$1<l<\infty$ and $F\in H^p$. As $s$ is continuous on
$\overline{\D}$, hence the existence of the non tangential
limit \eqref{limntG} and the majorization \eqref{maxinegw}
follow from the corresponding properties of $H^p$ functions
\cite[thm 3.1]{gar}.
From Fatou's lemma, we then get the first half of (\ref{Fatou}),
and since $\|F\|_{L^p(\T_r)}\leq\|\tr\, F\|_{L^p(\T)}$ for 
$F\in H^p$ 
\cite[thm 1.5]{duren}, we obtain by \eqref{estims}
$$
\|w\|_{L^p(\T_r)}\leq e^{2\|s\|_{L^\infty( \D)}}\|\tr\, w\|_{L^p(\T)}
\leq e^{8\|\alpha\|_{L^\infty(\D)}}\|\tr\, w\|_{L^p(\T)},
$$
which yields the second half of (\ref{Fatou}). Finally, (\ref{limgp}) follows 
from the continuity of $s$ and the corresponding property for $H^p$-functions
\cite[thm. 2.6]{duren}. This demonstrates assertion 1.

Since $e^s$ is continuous and never zero on $\D$, as noticed in 
Remark \ref{rems}, assertion 2 is a consequence of \eqref{factBN} and
of the corresponding properties for $H^p$-functions 
\cite[thms. 2.2, 2.3]{duren}. 

We turn to the proof of assertion $3$. Assume first that $w\in L^p(\D)$
satisfies $w=\cC \varphi+T(\alpha\overline{w})$ for some $\varphi\in L^p(\T)$.
As $\bar\d\cC\varphi=0$ on $\D$ because $\cC\varphi$ is holomorphic there,
we know from Proposition \ref{properties},  point 4,   that 
$\bar\d w=\alpha\bar w$ on $\D$. Further, the M. Riesz theorem yields
$$
\|\cC\varphi\|_{H^p(\D)}\leq C_p\|\varphi\|_{L^p(\T)}
$$
hence $w-T(\alpha\overline{w})\in H^p$. Lemma \ref{improve} now
provides us with the chain of inequalities:
$$
\|T(\alpha\overline{w})\|_{L^p(\T_r)}\leq \|T(\alpha\overline{w})\|_{L^\infty(\D)}\leq C'_p
\|T(\alpha\overline{w})\|_{W^{1,p^*}(\D)}\leq C_p''\left(\|w\|_{L^p(\D)}+\|\cC\varphi\|_{H^p(\D)}\right),
$$
where we used the Sobolev imbedding theorem. 
Therefore $w\in G_\alpha^p$ and (\ref{gpnorm}) holds.

Conversely, let $w\in G^p_{\alpha}$. 
Then $\bar\d w=\alpha\bar w\in L^p(\D)$ and Proposition \ref{properties}, 
assertion 4, tells us that the $L^p(\D)$-function $w-T(\alpha \bar w)$
anihilates $\bar\d$
in the distributional sense, hence is holomorphic on $\D$ by
Weyl's lemma. From  Proposition \ref{properties}, point 4 again,  
this entails $w\in W^{1,p}_{loc}(\D)$.
Appealing to Lemma \ref{gc},
with $r\D$ in place of $\D$, we obtain
$$
w(z)={1\over 2\pi i}\int_{\T_r}\frac{w(\zeta)}{\zeta-z}d\zeta+
T(\alpha\overline{w}\chi_{\D_r}) (z)\,,\ \ 
 |z|<r<1,
$$
and letting $r\to 1$ we get by  dominated convergence, (\ref{limgp}), 
and Proposition
\ref{properties}, point 4, that
\begin{equation}\label{Fatou'}
w=\cC(\tr\, w)+T(\alpha\overline{w}) \quad \mbox{ a. e. in }\D 
\, .
\end{equation}
Assertion 3 is now completely proven. 

Finally, assume that 
$w\in G^p_{\alpha}$ satisfies $\mbox{Re} \,\tr\, w=0$.
By Theorem \ref{trick2} we can write $w=e^sF$,
where $s$ is continuous on $\overline{\D}$ and real on $\T$, while
$F\in H^p(\D)$. Thus  $\mbox{Re} \,\tr\, F$ is zero on $\T$, and by the Poisson 
representation of $H^p$-functions  it follows that
$\mbox{Re}\, F\equiv0$ on $\D$ thus
$F$ is a pure imaginary constant, say, $c$ \cite[Thm 3.1]{duren}. 
Since $w = c \, e^s$, it has zero mean on $\T$ if and only if $c=0$.
When \eqref{alphabeta} holds, condition \eqref{normal:w}
likewise implies that $c=0$. \hfill\fin\par

\medskip

\subsection{Comparison between Sobolev and Hardy solutions} \label{proofpropf}

In this section, we prove Proposition \ref{hardysobol}. \par

\medskip

{\it Proof of Proposition \ref{hardysobol}.} \Rd We first check that $W^{1,p}(\D)\subset H^p_{\nu}(\D)$ and that (\ref{inegHSobt}) holds. \Bk Define $\alpha$ through 
\eqref{alphabeta}. In view of \eqref{correspf} and
Proposition  \ref{trick-hardy}, point (b), 
it is enough to check the corresponding 
property for $G^p_\alpha$. But if $w\in W^{1,p}(\D)$ meets (\ref{eq:w}), 
then Lemma \ref{gc} yields 
$w={\mathcal C}(\mbox{tr }w)+T(\alpha\overline{w})$, 
and since $\tr\, w\in W^{1-1/p,p}(\T)\subset L^p(\T)$ we get from Proposition
\ref{properties}, point 1, that 
$w-T(\alpha\overline{w})={\mathcal C}(\mbox{tr }w)$ lies in $H^p$, 
implying that $w\in G^p_\alpha$ by Remark \ref{rmkY}. Moreover,
it follows from \eqref{gpnorm} and the trace theorem that
$\|w\|_{G^p_\alpha}\leq C_{p,\alpha} \|w\|_{W^{1,p}(\D)}$, as desired.\par

\medskip

\noindent \Rd To end the proof of Proposition \ref{hardysobol}, we establish the more general fact that, if $f\in W^{1,p}(\D)$ has a nontangential limit almost everywhere on $\T$, then this nontangential limit coincides with the trace of $f$ in the Sobolev sense. We provide an argument because we could not 
locate this ``elementary'' result in the literature. Note that, \Bk when $p>2$, each $f\in W^{1,p}(\Omega)$ (has a representative which) 
extends continuouly to $\overline{\Omega}$ by the Sobolev imbedding theorem,
so the nontangential limit exists everywhere 
on $\partial \Omega$ and is in fact an unrestricted limit.\par
\noindent \Rd Assume now that $p\leq 2$. Notice that \Bk $f$ is the restriction to
$\D$ of some $\widetilde{f}\in W^{1,p}(\RR^2)$ with compact support
\cite[{ Thm} 4.26]{Adams}. By H\"older's inequality 
$\widetilde{f}\in W^{1,s}(\R^2)$ with 
some $s\in(1,2)$,  hence the non-Lebesgue points of $\widetilde{f}$
have Hausdorff 1-dimension zero \cite[{ Thms 3.3.3, 2.6.16}]{Ziemer}.
In particular, for \emph{each} $r>0$ we have that $re^{i\theta}$ is a
 Lebesgue point of $\widetilde{f}$ for a.e. $\theta$. Then, 
regularizing $\widetilde{f}$ yields a sequence
of $\cD(\R^2)$-functions converging to $\widetilde{f}$ both in
$W^{1,p}(\R^2)$ and pointwise a.e. on every circle $\T_r$
\cite[{ Thm} 1.6.1]{Ziemer}. Consequently 
if we put $f_r(\xi):=f(r\xi)$ for $0<r<1$ and $\xi\in\D$,
we deduce that $\tr\,f_r(\zeta)=f(r\zeta)$ for a.e.
$\zeta\in\T$, and since $f_r$ tends to $f$ in $W^{1,p}(\D)$ 
as $r\to1$ we get that $f(r\zeta)\to\tr\,f(\zeta)$ in
$W^{1-1/p,p}(\T)$. Finally,  
as any $W^{1-1/p,p}(\T)$-converging sequence
has a pointwise a.e. converging subsequence,
there is $r_n\to1$ such that $f(r_n\zeta)\to \tr\,f(\zeta)$ 
for  a.e. $\zeta\in\T$, hence
$\tr\,f$ must be the radial ({\it a fortiori} nontangential) 
limit of $f$ when the latter exists.

\Gr \begin{remarquesubsect}
Since Sobolev functions can be redefined on a set of zero measure so as
to be absolutely continuous on almost every line (\cite[Rem. 2.1.5]{Ziemer}), it is easy to see
(use polar coordinates) that any function in $W^{1,p}(\D)$ has a radial limit
almost everywhere on $\T$. Note, however, that a function in  $W^{1,2}(\D)$ may
have no nontangential limit at all (see \cite{carleson}). Note also that a more precise result involving capacity was proven for continuous functions in $W^{1,p}(\D)$, see \cite{m,r}.
\medskip

\end{remarquesubsect}\Bk

\subsection{The Dirichlet problem in the class  $G^p_{\alpha}$} \label{proofdirw}
{\it  Proof of Theorem \ref{dirichletp+}.} 
Let $g\in H^p$. By assertion 3 of Proposition \ref{gc-hardy} and 
Remark \ref{rmkY}, together with the 
Cauchy formula for $H^p$-functions, we see that $w\in L^p(\D)$ belongs to  
$G_\alpha^p$ and satisfies $P_+(\tr\, w)=\tr\, g$ if, and only if
$w-T(\alpha\overline{w})=g$. But since $g$ {\it a fortiori} lies in
$L^p(\D)$, there is a unique $w\in L^p(\D)$ to meet the latter equation
as follows from Proposition \ref{properties}, assertion 4.
Moreover, by the same assertion, it holds that
\[\|w\|_{L^p(\D)}\leq C_{p,\alpha}\|g\|_{L^p(\D)}\leq 
C_{p,\alpha}\|g\|_{H^p},\]
hence \eqref{estimwphi} holds in view of \eqref{gpnorm}.~\hfill\fin

\medskip

{\it Proof of Theorem \ref{thm:dirichletw}.}  For each pair
$(\varphi,c)\in L^p_\R(\T)\times\R$, set 
\[A(\varphi,c):=\left(\mbox{Re}\,\bigl(\tr\, w_{\varphi.c}\bigr)\,\,,\,
\mbox{Im}\frac{1}{2\pi}\int_0^{2\pi} 
\tr\, w_{\varphi.c}(e^{i\theta})\,d\theta\right)\in L^p_\R(\T)\times\R\,,
\]
 where $w_{\varphi,c}$ 
is the unique function
in $G_\alpha^p$ such that $P_+( \tr\, w_{\varphi,c})=\varphi+i(\cH_0\varphi+c)$.

Observe from the M. Riesz theorem that
$\varphi\mapsto\varphi+i\cH_0\varphi$ is continuous from 
$L^p_{\R}(\T)$ into $\tr\, H^{p,0}\subset L^p(\T)$,
hence $A$ is well-defined and continuous 
from $L^p_\R(\T)\times\R$ into itself by \eqref{Fatou} and 
Theorem \ref{dirichletp+}. 

Put for simplicity 
$T_\alpha(w)=T(\alpha\overline{w})$. In view of \eqref{formula}, we have that
$(I-T_\alpha)w_{\varphi,c}=g$
where $g\in H^p$ satisfies $\tr\, g=\varphi+i(\cH_0\varphi+c)$, hence
$A(\varphi,c)=(\varphi,c)+B(\varphi,c)$ where
\[B(\varphi,c):=\left(
\mbox{Re}\,\bigl(\tr\, T_\alpha( w_{\varphi,c} )\bigr)\,\,,
\,\, 
\mbox{Im}\frac{1}{2\pi}\int_0^{2\pi} 
\tr\,  w_{\varphi,c}  (e^{i\theta})\,d\theta\,-c\right).\]
Since the $G_\alpha^p$-norm is finer than the $L^p(\D)$-norm,
we see from Theorem \ref{dirichletp+} and
Proposition \ref{properties}, point 4, that
$(\varphi,c)\mapsto T_\alpha(w_{\varphi,c})$ is continuous
from $L^p(\T)\times\R$ into $W^{1,p}(\D)$, so 
the first component of $B$ is continuous from 
$L^p(\T)\times\R$ into $W^{1-1/p,p}(\T)$. Therefore it
is compact from $L^p(\T)\times\R$ into $L^p(\T)$ and, since the 
second component is $\R$-valued and continuous, $B$ is compact
from $L^p(\T)\times\R$ into itself. Moreover, $A$ is injective by
Proposition \ref{gc-hardy}, point 4, consequently it is an isomorphism of
$L^p(\T)\times \R$ \Gr (see {\it e.g.} 
\cite[Ch. XVII, prop. 2.3]{Torchinsky}) \Bk thereby establishing the existence and uniqueness
part of (a). 
To establish \eqref{inegDw},
put $\tr\, w=\psi+iv$, where $v\in L^p_\R(\T)$ in view of 
Proposition \ref{gc-hardy}. Thanks to Theorem \ref{trick2},
we can write $\psi+iv=e^{{\rm tr }\,s}F$ where $F\in H^p(\D)$ and 
$\tr\,s$ is real-valued, $s\in C^{0,\gamma}(\overline{\D})$ for 
$0\leq \gamma<1$ ({\it cf.} Remark
\ref{rems}).
By definition $\tr\, F=h+i {\cH}_0 h +ib$, where 
$h\in L_\R^p(\T)$
and $b$ is a real constant. Thus $v=e^{{\rm tr}\, s}({\cH}_0 h+b)$ and 
$\psi=e^{{\rm tr}\, s}h$, which gives us
\begin{equation}
\label{vexpH}
v=e^{{\rm tr}\, s} {\cH}_0 (e^{-{\rm tr}\, s}\psi) +be^{{\rm tr}\, s}.
\end{equation}
Since, $v$ has mean $c$ on $\T$,
we get from \eqref{estims}
and the M. Riesz theorem that
$$
|b|\leq e^{4\|\alpha\|_{L^\infty(\D)}}\left(
c_pe^{8\|\alpha\|_{L^\infty(\D)}}
\|\psi\|_{L^p(\T)}+ |c|\right),
$$
where $c_p$ is the norm of ${\cH}_0$ on $L^p(\T)$.
Plugging this in \eqref{vexpH} implies now
$$
\|v\|_{L^p(\T)}\leq e^{8\|\alpha\|_{L^\infty(\D)}}
\Bigl(|c|+\, \left(1+e^{8\|\alpha\|_{L^\infty(\D)}}\right)c_p
\|\psi\|_{L^p(\T)}\Bigr)
$$
from which \eqref{inegDw} follows immediately. This concludes
the proof of (a).

Assume next that \eqref{alphabeta} holds and let 
$w_1\in G^p_\alpha$, satisfy $\mbox{Re}(\tr\, w_1)=\psi$ a.e. on $\T$. Such a 
$w_1$ exists by (a). From the observation made after Proposition 
\ref{trick-hardy}, we see that
the function $w$ defined as
\[w:=w_1-i\sigma^{-1/2}
\frac{1}{2\pi}\int_0^{2\pi} \left(\sigma^{1/2} \mbox{Im} \, {\tr\, w_1 }\right)(e^{i \theta}) \, d   
\theta\]
lies in $G^p_\alpha$,  and since it readily
satisfies \eqref{normal:w} we deduce that $w\in G^{p,0}_\alpha$.
Clearly $w$ has the same real part as $w_1$,
and by Proposition \ref{gc-hardy} point 4 it is the only member of 
$G^{p,0}_\alpha$ with this property.
This settles the existence and uniqueness part in (b).

The reasoning leading to
\eqref{inegDw} is easily adapted to yield \eqref{inegDw0},
upon trading the mean-equal-to-$c$ condition for \eqref{normal:w}
and taking \eqref{ellipticsigma} into account. This completes the proof.  \hfill\fin

\subsection{\Gr A Fatou theorem for $\mbox{Re }H^p_{\nu}$ \Bk} \label{proofdirf}

As a preparation for the proof of Theorem \ref{thm:dirichletu}, we  
establish the 
following Hodge type lemma.
\begin{lemmasubsect}
\label{Hodgelisse}
Define two subspaces ${\mathcal G}_{0,\infty}$ and
${\mathcal D}_{0,\infty}$
of $C^\infty_\R(\overline{\D})\times C^\infty_\R(\overline{\D})$ by
\[{\mathcal G}_{0,\infty}:=
\{\nabla g;\ g\in C^\infty_\R(\overline{\D}),\ \tr\, g=0\},\ 
{\mathcal D}_{0,\infty}:=\{(\partial_y h,-\partial_x h)^T;\ h\in 
C^\infty_\R(\overline{\D}),\ \tr\, h=0\}.
\]
Then, each $V\in C^\infty_\R(\overline{\D})\times C^\infty_\R(\overline{\D})$
that vanishes on $\T$ can be written 
uniquely in the form $V=G+D$, 
where $G\in {\mathcal G}_{0,\infty}$ and
$D\in {\mathcal D}_{0,\infty}$. Moreover, it holds 
for some constant $C_p$ that
\begin{equation}
\label{inegHodgep}
\|G\|_{L^p(\D)}+\|D\|_{L^p(\D)}\leq  C_p\|V\|_{L^p(\D)},
\end{equation}
where the subscript $L^p(\D)$ refers here to the norm of an $\R^2$-valued 
mapping.
\end{lemmasubsect}
{\it Proof.}  Recall from \eqref{defHodge} the Hodge decomposition
$V=G+D$, with
$G\in {\mathcal G}_{0,p}$ and $D\in {\mathcal D}_p$, for which
\eqref{inegHodgep} is known to hold \cite[{ Thm} 10.5.1]{im2}. 
Put $G=(\partial_x g,\partial_y g)^T$ and 
$D=(\partial_y h,-\partial_x h)^T$ with $g\in W^{1,p}_{0,\R}(\D)$ and
$h\in W^{1,p}_\R(\D)$.
Since $V\in C^\infty_\R(\overline{\D})\times C^\infty_\R(\overline{\D})$, 
the same is true of $G$ and $D$ \cite[Sec. 10.5]{im2}, hence
$g,h\in C^\infty_\R(\overline{\D})$ and $\tr\, g=0$. Then, all we have to show 
is that $h_{|_\T}$ is constant, because substracting this constant will 
produce a new $h$ with vanishing trace on $\T$, as desired.

Now, because $\tr\, g=0$ we deduce that $G$ is normal to $\T$ there,
and since $V_{|_\T}=0$ it follows that $D$ is also a normal vector field
on $\T$. Consequently $x\partial_y h(x,y)-y\partial_x h(x,y)=0$ 
for $x+iy\in\T$, which means exactly that $h$ is constant on $\T$.
\hfill\fin\par

{\it Proof of Theorem \ref{thm:dirichletu}.}
Since $\partial_x (\sigma\partial_x u)= 
\partial_y (-\sigma\partial_y u)$ by  \eqref{div},  
there is a distribution $v$ 
on $\D$ such that \eqref{system} holds  
\cite[Ch. II, Sec. 6, { Thm} VI]{Schwartz}. Then, for $\Phi\in{\mathcal D}_\R(\D)$,
we obtain
\[
\begin{array}{l}
\langle v,\partial_x\Phi\rangle=\langle\sigma\partial_yu,\Phi\rangle=-\langle u,\sigma \partial_y\Phi+
\Phi\partial_y\sigma\rangle,\\
\langle v,\partial_y\Phi\rangle=-\langle\sigma\partial_xu,\Phi\rangle=\langle u,\sigma \partial_x\Phi+
\Phi\partial_x\sigma\rangle,
\end{array}
\]
which entails by \eqref{ellipticsigma}, H\"older's inequality, and the 
Poincar\'e inequality that
\begin{equation}
\label{inegdLp}
\begin{array}{l}
|\langle v,\partial_x\Phi\rangle|\leq \|u\|_{L^p(\D)} C_{p,\sigma}\|\nabla\Phi\|_{L^q(\D)} \, , \ 
|\langle v,\partial_y\Phi\rangle|\leq \|u\|_{L^p(\D)} C_{p,\sigma}\|\nabla\Phi\|_{L^q(\D)}.
\end{array}
\end{equation}
Next, we observe that any $g\in C^\infty_\R(\overline{\D})$ 
satisfying $\tr\, g=0$ lies in $W^{1,q}_{0,\R}(\D)$, therefore it is the
limit in $W^{1,q}_\R(\D)$ of some sequence $\Phi_n\in\cD_\R(\D)$. In 
particular $\nabla\Phi_n$ converges to $\nabla g$ in $L^q(\D)$, 
implying in view of \eqref{inegdLp} that 
$(\Phi_1,\Phi_2)\mapsto 
\langle v,\Phi_1+\Phi_2\rangle$  is a bounded functional on 
both ${\mathcal G}_{0,\infty}$ and
${\mathcal D}_{0,\infty}$ when endowed with the $L^q(\D)$-norm.
By Lemma \ref{Hodgelisse}, this functional is \Gr $L^q(\D)\times L^q(\D)$ \Bk bounded on
the subspace of $C^\infty_\R(\overline{\D})\times C^\infty_\R(\overline{\D})$
comprising those vector fields that vanish on $\T$. Therefore,
by density, $(\Phi_1,\Phi_2)\mapsto 
\langle v,\Phi_1+\Phi_2\rangle$  is a bounded functional on 
$L^q_\R(\D)\times L^q_\R(\D)$, so that in fact $v\in L^p_\R(\D)$.
If we put $f=u+iv\in L^p(\D)$ and 
$\nu=(1-\sigma)/(1+\sigma)\in W^{1,\infty}_\R(\D)$, it is now a mechanical 
consequence of 
\eqref{system} that equation \eqref{dbar} is satisfied in the distributional 
sense. Defining $w\in L^p(\D)$ and $\alpha\in L^\infty(\D)$ 
through \eqref{correspw} and \eqref{alphabeta}, we see from Proposition
\ref{trick-hardy} that $w$ solves \eqref{eq:w}. 
Hence Theorem \ref{trick2} applies to the effect that $w$, thus also $f$
and {\it a fortiori} $u$, lie
in $W^{1,l}_{loc}$ for $l\in(1,\infty)$. Moreover, 
we may  write \eqref{factBN} with $\mbox{Im}\,\tr\, s=0$ and, say,
$F=a+ib$ a holomorphic function in $\D$. As $s$ lies in 
$C^{0,\gamma}(\overline{\D})$, for  each $\varepsilon>0$ we can pick 
$r_0$ such that
$|\mbox{Im} \, \exp(s(z))|<\varepsilon |\exp (s(z))|$ as soon as $r_0<|z|\leq1$,
and for such $z$ we deduce from \eqref{estims} that
\[\mbox{Re} \, w(z)\geq e^{-4\|\alpha\|_{L^\infty(\D)}}((1-\varepsilon^2)^{1/2}|a(z)|-\varepsilon|b(z)|).\]
Since $b_r=\cH_0(a_r)$ ( recall that $b_r(z)=b(rz)$, resp. $a_r(z)=a(rz)$, for all $z\in \T$), we obtain from the M. Riesz theorem
when $r_0<r<1$ that
\[\|\mbox{Re}\, w(z)\|_{L^p(\T_r)} \geq e^{-4\|\alpha\|_{L^\infty(\D)}}
((1-\varepsilon^2)^{1/2}-\varepsilon c_p)\|a(z)\|_{L^p(\T_r)},\]
and picking $\varepsilon$ so small that
$((1-\varepsilon^2)^{1/2}-\varepsilon c_p)=C>0$ we conclude 
by \eqref{correspw} that
\[\|b\|_{L^p_r(\T)}\leq c_p\|a\|_{L^p_r(\T)}\leq
C_{p,\sigma,u}\|\mbox{Re} \, w(z)\|_{L^p(\T_r)}\leq
C'_{p,\sigma,u}\|u\|_{L^p(\T_r)},\ \ \  r_0<r<1.\] 
Consequently, if $\|u\|_{F,p}<\infty$, then $F=a+ib\in H^p$ so that
$w\in G^p_\alpha$ by Theorem \ref{trick2}, whence $f\in H^p_\nu$ by
Proposition \ref{trick-hardy}. All the assertions now readily follow from
Theorems \ref{gc-hardy} and \ref{thm:dirichletf}.
\hfill\fin\par

\subsection{Higher regularity} \label{proofhigherreg}

In order to prove Theorem \ref{sobnontang}, we shall
make use of the following observation:
\begin{lemmasubsect}
\label{Beltcd}
Let $\varphi\in W_{\R}^{1-1/p,p}(\T)$, and $f\in W^{1,p}(\D)$ 
be the unique solution to \eqref{dbar} satisfying
${\rm Re} (\tr\, f)=\varphi$ and \eqref{normal:f}, 
{\it cf.} Theorem \ref{dirichlet}. 
Then $W=(1-\nu^2)^{1/2} \, \partial f\in L^p(\D)$ is a solution 
to \eqref{eq:w} with $\alpha = \partial\nu/(1-\nu^2)$, 
that is to say
\begin{equation}
\label{eqbder}
\overline{\partial} W = \frac{\partial\nu}{1-\nu^2}\overline{W}
\end{equation}
in the sense of distributions.
Moreover, it holds that
\begin{equation}
\label{BNRt}
\partial f= e^{s} F
\end{equation}
where $F$ is holomorphic in $\D$,
$s\in C^{0,\gamma}(\overline{\D})$ for 
every $0<\gamma<1$, and for some constant $c_\kappa$ we have
$\|s\|_{L^{\infty}(\D)}\leq c_\kappa
\|\nu\|_{W^{1,\infty}(\D)}$.
\end{lemmasubsect}
{\it Proof.} As $\nu\in W^{1,\infty}(\D)$, observe that the distributional derivative of 
$\nu\overline{\partial f}\in L^p(\D)$ can be computed 
according to Leibniz's rule: 
\begin{equation}
\label{Leibniz}
\partial\left(\nu\overline{\d f}\right) =
\partial \nu \overline{\d f}+\nu\partial\left(\overline \d \overline{f}\right),
\end{equation}
where we emphasize that the second summand in the right-hand side 
of \eqref{Leibniz} is 
to be interpreted as indicated in the footnote before Theorem \ref{dirichlet}. \Bk Indeed, pick a function $\varphi\in {\mathcal D}(\D)$. By definition,
$$
\langle \nu\overline{\partial f},\partial\varphi\rangle=-\langle \nu\overline{f},\overline{\partial}\partial\varphi\rangle-\langle \overline{\partial}\nu \overline{f},\partial\varphi\rangle.
$$
By the Leibniz rule
$$
-\langle \nu\overline{f},\overline{\partial}\partial \varphi\rangle  =  -\langle \overline{f},\nu\overline{\partial}\partial\varphi\rangle
 =  -\langle \overline{f},\overline{\partial}(\nu\partial\varphi)-\overline{\partial}\nu\partial\varphi\rangle.
$$
It follows that
$$
\begin{array}{l}
\langle \nu\overline{\partial f},\partial\varphi\rangle  =  -\langle \overline{f},\overline{\partial}(\nu\partial\varphi)\rangle =  \langle \overline{\partial}\overline{f},\nu\partial\varphi\rangle  =  \langle \overline{\partial}\overline{f},\partial(\nu\varphi)-\varphi\partial\nu\rangle\\
 =  -\langle \partial\overline{\partial}\overline{f},\nu\varphi\rangle-\langle \partial\nu\overline{\partial}\overline{f},\varphi\rangle  =  -\langle \nu\partial\overline{\partial}\overline{f},\varphi\rangle-\langle \partial\nu\overline{\partial}\overline{f},\varphi\rangle,
\end{array}
$$
which is the desired conclusion. \Bk\par
\noindent Setting $G:=\partial f$ and applying $\partial$ to \eqref{dbar}, 
we thus obtain, since $\partial$ and
$\overline \partial$ commute, that
$\overline{\partial} G=\nu \partial\overline{G}+
(\partial\nu)\overline{G}$.
As $\nu$ is real, conjugating this last equation provides us with
an expression for $\partial\overline{G}$,
and solving for $\overline{\d} G$ after substituting back yields 
\[
 \overline{\partial} G=\frac{\nu\overline{\partial}\nu}{1-\nu^2}G+
 \frac{\partial\nu}{1-\nu^2}\overline{G} \, 
\]
from which we deduce that $W = (1-\nu^2)^{1/2} \, G$ satisfies
(\ref{eqbder}), as $\overline{\d} W$ can in turn be computed by the chain rule 
because $(1-\nu^2)^{1/2}\in W^{1,\infty}(\D)$.
The remaining assertions follow from Theorem \ref{trick2}
upon setting  $\alpha=\partial\nu/(1-\nu^2)\in L^\infty(\D)$.
\hfill\fin
\par
\medskip

\par
\medskip
\noindent {\it Proof of Proposition \ref{hnusob}.} Let 
$u\in W^{1,p}_\R(\T)\subset W_\R^{1-1/p,p}(\T)$ and 
put $v={\mathcal H}_{\nu}u$. We must show that 
$v\in W^{1,p}_\R(\T)$, and for this we may assume that $u$ has zero mean 
on $\T$ for adding a constant to $u$ does not affect $v$. Then,
$v={\mathcal H}_{\nu}u$ becomes equivalent to
$u=-{\mathcal H}_{-\nu}v$ as follows immediately from the fact that
$f$ satisfies \eqref{dbar} if, and only if $if$ satisfies a similar
equation with $\nu$ replaced by $-\nu$.
From Theorem \ref{dirichlet} we know that
$v\in W_\R^{1-1/p,p}(\T)$, so let $v_n$ be a 
sequence of $C^\infty(\T)$-functions converging to $v$ there and 
set $u_n=-{\mathcal H}_{-\nu}v_n$. 
Since $v_n$ converges to $v$ in $W^{1-1/p,p}(\T)$, we 
get from Corollary \ref{MRieszgen} that $u_n$ converges to $u$ there.
By the definition of ${\mathcal H}_\nu$, we have
that $u_n+iv_n$ is the trace on $\T$ of the solution to 
\eqref{dbar} in $W^{1,p}(\D)$ whose real part on $\T$ is $u_n$.
With a slight abuse of notation, we still designate by $u_n$, $v_n$ 
the real and imaginary parts of that solution in $W_\R^{1,p}(\D)$. 
By inspection, the generalized Cauchy-Riemann equations
\eqref{system} do hold with $u$ replaced by $u_n$ and $v$ by $v_n$,
so that $u_n$, $v_n$ may as well be characterized respectively
as the unique solutions in $\D$  to \eqref{div}, \eqref{system2} 
whose traces on $\T$ are our previous $u_n$ and $v_n$
\cite{campanato}.

Now, we know that $u_n$, $v_n$ lie in
$W^{1,p}(\D)$; however, since $u_n$ is smooth on $\T$,
it follows from \cite[Thm 9.15]{gt} that in fact
$u_n\in W_\R^{2,r}(\D)$  for all $r\in (1,\infty)$. In view of \eqref{system},
we deduce that the same is true of $v_n$ as 
$\sigma\in W_\R^{1,\infty}(\D)$.
Pick any $\varphi\in C^\infty(\T)$ and put $\psi={\mathcal H}_\nu\varphi$
with $\nu$-conjugate $\psi$. The $W^{2,r}(\D)$ regularity just mentioned
allows us by density to apply the divergence formula 
so as to obtain 
\[\int_{\T} v_n\partial_{\theta}\varphi\,d\theta=
-2\pi \int_{\T} v_n\partial_{n}\psi/\sigma\,d\theta=
-2\pi \int_{\T} \partial_n v_n\psi/\sigma\,d\theta=
\int_{\T} \partial_\theta u_n\psi\,d\theta,
\]
where we used \eqref{CRtang}.
As  $u_n$, $v_n$ converge to $u$, $v$ in $W^{1-1/p,p}(\T)$,
we get in the limit, since differentiation is continuous 
$W^{1-1/p,p}(\T)\to W^{-1/p,p}(\T)$, that
\[\int_{\T} v\partial_{\theta}\varphi\,d\theta=
\int_{\T} \partial_\theta u\psi\,d\theta,
\]
where $\partial_\theta u$ is to be understood as a member of $W^{-1/p,p}$.
However, we have by assumption that in fact 
$\partial_\theta u\in L^p(\T)$,
therefore from H\"older's inequality
\[
\left|\int_{\T}v\partial_\theta \varphi\,d\theta\right|\leq 
\|\partial_\theta u\|_{L^p(\T)}\|\psi\|_{L^q(\T)}.
\]
But from Corollary \ref{MRieszgen} we know that
$\|\psi\|_{L^q(\T)}\leq C_\nu \|\varphi\|_{L^q(\T)}$,
so the distribution $\partial_{\theta}v$ 
in fact lies in $L^p(\T)$ and the conclusion holds. \hfill\fin\par
\medskip
\medskip

{\it Proof of Theorem \ref{sobnontang}.} Write $f=u+iv$ to indicate the real and imaginary parts of $f$. By Proposition \ref{hnusob}, $\tr\, f\in W^{1,p}(\T)$. First we shall prove that the left
hand sides of (\ref{HCp1}) and (\ref{HCp2}) are finite, that is,
we show $\partial f$
and $\overline{\partial} f$ satisfy a Hardy condition of order $p$. 
To this effect, we consider $w=(f-\nu \bar{f})/\sqrt{1-\nu^2}$
and we establish the equivalent fact 
that both  $\partial w$ and $\bar{\partial} w$ satisfy a
Hardy condition of order $p$.  We first notice:
\begin{lemmasubsect} \label{improvew}
If $p_1\in[p,2p)$ then $f,w\in W^{1,p_1}(\D)$;
in particular $f,w \in C^{0,1- 2/p_0}(\overline{\D})$
for some $p_0>2$.
\end{lemmasubsect}
{\it Proof.} 
Note that $\tr w\in W^{1,p}(\T)$ since $\tr f$ does.
\Rd As we will see in the proof of  Proposition \ref{properties} in Appendix \ref{technic} ({\it cf.} 
(\ref{boundC2})), \Bk it entails  
$\partial (\mathcal{C}\tr w) \in H^p$, and since $\overline{\partial} (\mathcal{C}\tr w)=0$
we see from the Poincar\'e inequality and Lemma \ref{improve}, point 1,
that  $(\mathcal{C}\tr w)$ lies in $W^{1,p_1}(\D)$. 
As $w$, thus $\alpha \overline{w}$ belongs to $L^{p_1}(\D)$ by Theorem
\ref{trick2}, we get from Proposition \ref{properties}, point 4,
that $T_\alpha w \in W^{1,p_1}(\D)$. Since
$$
w={\mathcal C}(\tr w)+T_{\alpha}w,
$$
we see that $w$, therefore also $f$ is in $w\in W^{1,p_1}(\D)$.
As $2p>2$, the last assertion 
now follows from the Sobolev imbedding theorem. \hfill\fin\par
\medskip
Back to proof of Theorem \ref{sobnontang} we observe that,
to prove the finiteness of the left hand sides in (\ref{HCp1}) and 
(\ref{HCp2}), we may as well add a 
real constant to $f$. Since the latter is (even H\"older)
continuous on $\overline{\D}$ by Lemma \ref{improvew},
we may thus assume that its real part is 
larger than a positive constant.
Then, the same is true of $w=(f-\nu\overline{f})/\sqrt{1-\nu^2}$,
say,  ${\rm Re} \, w(z)\geq c_0>0$ for $z\in\D$.
This results in  $\overline{w}/w$ being  H\"older continuous of exponent
$1-2/p_0$
in $\D$. Now, consider the function $s$ introduced in Theorem \ref{trick2}.
Letting $B(z,\varepsilon)$ indicate the ball of center $z$ with radius 
$\varepsilon$, we gather from \eqref{defs} that, for a.e. $z\in\D$,
\begin{equation}
\label{calculps}
\partial s(z)=
\lim_{\varepsilon\to0}\frac1{2\pi i}\iint_{\D\setminus B(z,\varepsilon)}
\frac{r(\zeta)}{(\zeta-z)^2}d\zeta\wedge d\bar\zeta \,+
\frac1{2\pi i}\iint_\D\frac{\overline{ r(\zeta)}}
{(1-\bar\zeta z)^2}
d\zeta\wedge d\bar\zeta\, , 
\end{equation}
where the function $r$ was defined as $r=\alpha \overline{w}/w$
and the existence of the limit  a.e. comes from the existence of
the Beurling transform as a singular integral operator of 
Calder\`on-Zygmund type.  To evaluate the first integral in 
\eqref{calculps}, we establish a lemma which is best stated in terms of the 
space $BMOA(\D)$, 
comprised of those $H^2$-functions whose trace on $\T$
has bounded mean oscillation, see {\it e.g.} \cite[p. 240]{duren}.
To us, the important fact will be that $BMOA(\D)\subset H^p$ for all 
$p<\infty$. 
\begin{lemmasubsect} \label{summand}
There exist a function $b\in L^{\infty}(\D)$ and a function 
$\psi\in BMOA(\D)$ such that, for a.e.  $z\in \D$,
\begin{equation}
\label{decbpsi}
\lim_{\varepsilon\to0}\frac1{2\pi i}\iint_{\D\setminus B(z,\varepsilon)}
\frac{r(\zeta)}{(\zeta-z)^2}d\zeta\wedge d\bar\zeta =b(z)+\psi(z).
\end{equation}
\end{lemmasubsect}
\medskip

\noindent {\it Proof.} We may rewrite the first integral in the right hand side of
(\ref{calculps}) as
\[\frac1{2\pi i}\iint_{\D\setminus B(z,\varepsilon)}
\frac{r(\zeta)}{(\zeta-z)^2}d\zeta\wedge d\bar\zeta \,=\,
\frac1{2\pi i}\iint_{\D\setminus B(z,\varepsilon)}
\frac{\alpha(\zeta)((\overline{w}/ w)(\zeta)-(\overline{w}/ w)(z))}
{(\zeta-z)^2}d\zeta\wedge d\bar\zeta\]
\[+
\frac{(\overline{w}/ w(z))}
{2\pi i}\iint_{\D\setminus B(z,\varepsilon)}
\frac{\alpha(\zeta)}{(\zeta- z)^2}
d\zeta\wedge d\bar\zeta \, .
\]
Since $\alpha$ is bounded and $\overline{w}/w$  is H\"older continuous 
of  order $1-2/p_0$, 
the first integral in the right hand side is majorized by
\[C_{\nu,w}\iint_\D |\zeta-z|^{-1-2/p_0}d\zeta\wedge d\bar\zeta \leq
C_{\nu,w}\iint_{|\xi|\leq2} |\xi|^{-1-2/p_0}dm(\xi)<+\infty.\]
As to the second integral, we put for simplicity
$\Phi:=\log(\sigma^{1/2})$ and we recall from
\eqref{alphabeta} that $\alpha=\overline{\partial}\Phi$, whence
by Stoke's theorem
\begin{equation}
\label{estPhi}
\frac{1}{2\pi i}\iint_{\D\setminus B(z,\varepsilon)}
\frac{\alpha(\zeta)}{(\zeta- z)^2}
d\zeta\wedge d\bar\zeta =-\frac{1}{2i\pi}\int_\T \frac{\Phi(\zeta)d\zeta}{(\zeta-z)^2}
+\frac{1}{2i\pi}\int_{\partial B(z,\varepsilon)}
\frac{\Phi(\zeta)d\zeta}{(\zeta-z)^2}.
\end{equation}
Since
\[\frac{1}{2i\pi}\int_{\partial B(z,\varepsilon)}
\frac{\Phi(\zeta)d\zeta}{(\zeta-z)^2}=
\frac{1}{2i\pi}\int_{\partial B(z,\varepsilon)}
\frac{(\Phi(\zeta)-\Phi(z))d\zeta}{(\zeta-z)^2},
\]
and $\Phi$ is Lipschitz continuous with constant, say, $K$
(because $\sigma\in W^{1,\infty}(\D)$ and in view of \eqref{ellipticsigma})
we get
\[\left|\frac{1}{2i\pi}\int_{\partial B(z,\varepsilon)}
\frac{(\Phi(\zeta)-\Phi(z))\,d\zeta}{(\zeta-z)^2}\right|\leq 
\frac{K\varepsilon}{2\pi}\int_{\partial B(z,\varepsilon)}
\frac{|d\zeta|}{|\zeta-z|^2}=K.
\]
Thus the second integral in the right hand side of (\ref{estPhi})
is uniformly bounded. As to the first, we observe by 
the Lipschitz character of $\Phi$ that $\Phi_{|_\T}$ 
is absolutely continuous on $\T$
with derivative $\varphi:=\partial_\theta \Phi(e^{i\theta})$
which is bounded in modulus by $K$ for a.e. $\theta$. Therefore,
integrating by parts, we get
\[\frac{1}{2i\pi}\int_\T \frac{\Phi(\zeta)d\zeta}{(\zeta-z)^2}
=\frac{1}{2i\pi}\int_\T \frac{\varphi(\zeta)d\zeta}{(\zeta-z)}
\]
which is the Cauchy integral of a bounded function and therefore
belongs {$BMOA(\D)$} \cite[Ch. VI, Cor. 2.5]{gar}.
\hfill\fin\par

\begin{lemmasubsect} \label{hardysigma}
The function $\partial s$ satisfies a Hardy condition of order $l$ for all 
$l\in (1,+\infty)$.
\end{lemmasubsect}
{\it Proof.} Let $\psi$ be as in Lemma \ref{summand} and
$a(z)$ be the holomorphic integral vanishing at $0$
of $\psi$, namely  $\partial a(z)=\psi(z)$ and $\bar\partial a(z)=0$
for $z\in\D$. If we set
\[B(z):=\frac1{2\pi i}\iint_\D
\frac{r(\zeta)}{\zeta-z}
d\zeta\wedge d\bar\zeta \, , \mbox{ for } z \in \D \,,
\]
we find by \eqref{decbpsi} and Proposition \ref{properties}, point 4,
that the bounded function $B-a$ 
has bounded partial derivatives 
$\overline{\partial}(B-a)(z)=r(z)$ and
$\partial(B-a)(z)=b(z)$, thus it lies in $W^{1,\infty}(\D)$ hence it
is Lipschitz continuous.
But the second summand in the right hand side of (\ref{calculps}) is,
by our very construction, the derivative of the holomorphic function 
\[H(z)=\frac1{2\pi i}\iint_\D
\frac{z\overline{r(\zeta)}}
{1-\bar\zeta z}
d\zeta\wedge d\bar\zeta \, , \quad z \in \D \, 
\]
that vanishes at 0 and whose real 
part on $\T$ is $\Upsilon=-{\rm Re} \, B_{|_\T}$, see the discussion after 
\eqref{defs}. Writing $-B=-(B-a)-a$, we have
that $\Upsilon=\Upsilon_1+ \Upsilon_2$ where $\Upsilon_1$ is 
Lipschitz continuous and thus absolutely continuous with bounded derivative 
on $\T$, while $\Upsilon_2=-{\rm Re} \, (\tr\, a)$.
Consequently $\tr\, (H+a)=\Upsilon_1+i{\cH}_0\left(\Upsilon_1\right)$ where  we recall that ${\cH}_0$ denotes the usual conjugation operator. Now, $\Upsilon_1$ is
{\it a fortiori} in $W^{1,l}(\T)$ for all
$1<l<\infty$, therefore the same is true of ${\cH}_0\left(\Upsilon_1\right)$
by Proposition \ref{hnusob} applied with $\nu=0$.
Consequently $(H+a)^\prime=H^\prime+\psi$ lies in $H^l$ for all
$1<l<\infty$, and finally so does $H^\prime$ since it is the case of
$\psi\in {BMOA(\D)}$. Altogether, considering separately the summands in 
the right-hand side of
(\ref{calculps}) and recalling Lemma \ref{summand}, we have proven that
$\partial s$ satisfies a Hardy condition of any order in $(1,+\infty)$. \hfill\fin

\medskip

Now, let us turn to the holomorphic function
$F\in H^p(\D)$ in the factorization 
$w=e^s F$ of Theorem \ref{trick2}. We claim:
\begin{lemmasubsect} \label{Fprimehardy}
The function $F^{\prime}\in H^p(\D)$.
\end{lemmasubsect}
{\it Proof.} Observe,
since ${\rm Re} \, {s}=0$ on $\T$, that
$|F|=|w|$ there. Moreover,  by Lemma \ref{improvew} and Remark \ref{rems},
$F=e^{-s}w$ is
(H\"older) continuous on $\overline{\D}$ and it does not vanish there 
by our assumption on $w$. Therefore $F$ can have no inner factor 
in its inner-outer decomposition 
\cite[Ch. II, Cor. 5.7, Thms 6.2, 6.3]{gar} thus it is an outer function:
\begin{equation}
\label{Fout}
F(z)=\xi_0\exp\left\{\frac{1}{2\pi}\int_0^{2\pi}\frac{e^{i\theta}+z}
{e^{i\theta}-z}\,\log |w(e^{i\theta})|\,d\theta \right\}
\end{equation}
with $\xi_0$ a unimodular constant. As $\tr\, w\in W^{1,p}(\T)$ 
is bounded in modulus from above and below by strictly positive constants,
$\log|w(e^{i\theta})|$ also lies in $W^{1,p}(\T)$. Hence, in view of
\eqref{boundC2} \Rd in Appendix \ref{technic} below \Bk, the 
derivative of the holomorphic function of $z$ defined by
\[\frac{1}{2\pi}\int_0^{2\pi}\frac{e^{i\theta}+z}
{e^{i\theta}-z}\,\log |w(e^{i\theta})|\,d\theta=
-\frac{1}{2\pi}\int_\T \log |w(e^{i\theta})|\,d\theta
+\frac{1}{i\pi}\int_\T
\frac{\log |w(\xi)|}
{\xi-z}\,d\xi
\]
lies in $H^p$, and by (\ref{Fout}) so does the derivative of $F$ since the 
latter is bounded. That is, we have proven that $F^\prime\in H^p$. \hfill\fin

\medskip

Now, since
\[\partial w=e^s \partial s\,F+e^s F^\prime,\]
we see by the boundedness of $F$, $s$, and the Hardy character of
$\partial s$, $F^\prime$ just established that $\partial w$ satisfies a 
Hardy condition of order $p$.
Besides, $\bar{\partial} w=\alpha \bar{w}$ is bounded,
being the product of a $L^\infty(\D)$ function and a H\"older continuous 
one.
Thus both $\partial w$ and $\bar{\partial} w$ satisfy a
Hardy condition of order $p$, and since $w$ is bounded 
it follows that $\partial f$ and $\overline{\partial} f$ also satisfy 
a  Hardy condition of order $p$, that is, the left-hand sides of 
(\ref{HCp1}) and (\ref{HCp2}) are finite, as announced.

From this, in view of Lemma \ref{Beltcd},
we deduce that
$W=(1-\nu)^{1/2}\partial f$ lies in $G_{\alpha_1}^p$ with
$\alpha_1=\partial\nu/(1-\nu^2)$. Clearly $\alpha_1$ lies in $L^\infty(\D)$, hence
by \eqref{kappa} and the relation  
$\overline{\partial} f=\nu\overline{\partial f}$ we deduce from
Proposition \ref{gc-hardy} that $\partial f$ and $\bar{\partial} f$ have nontangential limits,
say, $\Phi_1$ and $\Phi_2=\nu_{|_\T}\overline{\Phi_1}$ 
to which $\partial f(r e^{i\theta})$ and
$\overline{\partial} f(r e^{i\theta})$ converge in $L^p(\T)$ 
as $r\to1$. It only remains for us to establish 
the explicit expression \eqref{relim} for $\tr \partial f$, because the latter 
readily implies that $\|\tr \partial f\|_{L^p(\T)}\leq C_{\nu}\|\tr f\|_{W^{1,p}(\T)}$,
hence assertions $(b)$ and $(c)$ of Theorem \ref{sobnontang} will follow from  
Proposition \ref{gc-hardy} as applied to $W$, and since we already pointed out that
$(a)$ is a rephrasing of Proposition \ref{hnusob} the proof will be complete.

To establish \eqref{relim},  observe by the absolute continuity on a.e. coordinate 
line characterizing Sobolev functions 
(choose polar coordinates) that
\[f(re^{i\theta_1})-f(re^{i\theta_2})=
-\int_{\theta_1}^{\theta_2}\left(r\sin\theta(\partial f+\overline{\partial} f)
-ir\cos\theta(\partial f-\overline{\partial} f)\right)(re^{i\theta})\,d\theta\]
for a.e. $r\in(0,1)$ and all $\theta_1,\theta_2$. Letting $r\to1$, we obtain
by the continuity of $f$ and the $L^p(\T)$-convergence
of $\partial f(r e^{i\theta})$ and
$\overline{\partial} f(r e^{i\theta})$ to $\Phi_1$ and $\Phi_2$ that
\[f(e^{i\theta_1})-f(e^{i\theta_2})=
-\int_{\theta_1}^{\theta_2}\left(\sin\theta(\Phi_1+\Phi_2)
-i\cos\theta(\Phi_1-\Phi_2)\right)(e^{i\theta})\,d\theta, 
\]
thereby showing that
\[\partial_{\theta}f(e^{i\theta})=ie^{i\theta}\Phi_1(e^{i\theta})
-ie^{-i\theta}\Phi_2(e^{i\theta})=
ie^{i\theta}\Phi_1(e^{i\theta})
-ie^{-i\theta}\nu (e^{i\theta})\overline{\Phi_1(e^{i\theta})}.
\]
Conjugating this identity, we obtain since $\nu$ is real-valued that
\[\partial_{\theta}\overline{f(e^{i\theta})}=-ie^{-i\theta}
\overline{\Phi_1(e^{i\theta})}+ie^{i\theta}
\nu (e^{i\theta})\Phi_1(e^{i\theta}),
\]
whence
\[\partial_{\theta}f(e^{i\theta})-
\nu(e^{i\theta})\partial_{\theta}\overline{f(e^{i\theta})}=
(1-\nu^2)ie^{i\theta}\Phi_1(e^{i\theta}).
\]
which is (\ref{relim}).

\hfill\fin

\bigskip

{\it Proof of Corollary \ref{fjr}. } 
Let $f=u+iv \in W^{1,p}(\D)$ be the solution to (\ref{dbar}) satisfying 
(\ref{normal:f}) granted by Theorem \ref{dirichlet}. By \eqref{formhardy}
and (\ref{system}),  $f$ meets
\begin{equation} \label{hardytype}
\supess_{0<r<1} \left\Vert \nabla f\right\Vert_{L^p(\T_r)}<+\infty.
\end{equation}
If $W:=(1-\nu^2)^{1/2}\partial f$, we saw in the proof of Lemma \ref{Beltcd} 
that  $\overline{\partial}W=\alpha\overline{W}$ for some 
$\alpha\in L^{\infty}(\D)$. Thus
$W\in G^p_{\alpha}(\D)$ by (\ref{hardytype}), hence it belongs both to
$W^{1,l}_{loc}(\D)$ for $1 < l < \infty$ and to
$L^{p_1}(\D)$, for $p_1 \in (p, 2p)$
in view of Theorem \ref{trick2}.
Clearly the same is true of
$\partial f$ and $\overline{\partial}f=\nu\overline{\partial f}$,
in particular $f \in W^{1,p_2}(\D)$ for some $p_2 > 2$, hence is H\"older continuous on $\overline{\D}$.
By the Sobolev embedding theorem as applied to its derivatives, $f$ is moreover continuously differentiable 
in $\D$, so we get for every $r e^{i\theta'} \in \D$ that:
\begin{equation}
\label{eqqq}
 f(r e^{i\theta'}) - f(r) = \displaystyle \int_0^{\theta'} (\partial_{\theta} f) (r e^{i\theta}) d \theta
\end{equation}
where $(\partial_{\theta} f) (z)  = i (z \, \partial f (z)- \bar{z} \nu \overline{\partial f}(z))$.
Now, by Proposition \ref{gc-hardy} applied to $W$, the derivatives
$\partial f$ and $\overline{\partial}f$ have
non tangential limits $h$ and $\nu \overline{h}$ respectively, 
where $h \in L^p(\T)$, and by (\ref{limgp}):
$$
\lim_{r\rightarrow 1} \int_0^{2\pi} \left\vert \partial f(re^{i\theta})-h(e^{i\theta})\right\vert^pd\theta=
\lim_{r\rightarrow 1} \int_0^{2\pi} \left\vert \overline{\partial} f(re^{i\theta})-\nu \, \bar{h}(e^{i\theta})\right\vert^pd\theta=0 \, .
$$
Passing to the limit as $r \to 1$ in (\ref{eqqq}) yields that
$\tr f$ is absolutely continuous on $\T$ with tangential derivative
$i (e^{i\theta} h - e^{-i\theta} \nu \bar{h}) \in L^p(\T)$, proving that
$\tr f\in W^{1,p}(\T)$. \par 
\noindent Since $\partial f\in G^p_{\alpha}(\D)$,
Proposition \ref{gc-hardy} implies ${\mathcal M}_{\partial f}\in L^p(\T)$. 
The same is true of 
${\mathcal M}_{\overline{\partial}f}$ because of (\ref{dbar}) and the boundedness of $\nu$, therefore ${\mathcal M}_{\left\Vert \nabla f\right\Vert}\in L^p(\T)$. Moreover, since $\partial f$ and $\overline{\partial f}$ both have non tangential limits in $L^p(\T)$, so does  $\nabla f$. The same conclusions 
then hold for $u=\mbox{Re}f$. In addition, the Green-Riemann formula ensures that
$$
0=\int_{\D} \mbox{div}(\sigma\nabla u)(x)dx=\int_{\T} \sigma\partial_n\tr u\,.
$$
This establishes point 1.

As to point 2,  let $v\in W^{1,p}_\R(\T)$ be such that $\partial_{\theta}v=\sigma g$. By Theorem \ref{dirichlet}, there exists $f\in W^{1,p}(\D)$ such that (\ref{dbar}) holds and $\mbox{Im }\tr f=v$ on $\T$. Then,
Theorem \ref{sobnontang} implies that $u:=\mbox{Re }f$ fulfills all the requirements. \par
\noindent To prove uniqueness up to an additive constant, consider $u\in W^{1,p}_{\R}(\D)$ with
$\mbox{div}(\sigma\nabla u)=0$ in $\D$ and $\partial_nu=0$ on $\T$. Put $f=u+iv$ for a solution
to \eqref{dbar} in $W^{1,p}(\D)$ such that $u=\mbox{Re }f$. By Remark \ref{remNeuman}
$\partial_\theta v=\sigma\partial_nu=0$ on $\T$, which means that $\tr v$ is constant,
{\it i.e.} $\mbox{Im }\tr f$ is constant.
By Theorem \ref{dirichlet}, $u$ is constant in $\D$. \hfill\fin \Bk
\subsection{Proofs of the density results} \label{proofdense}
\subsubsection{Density for Sobolev traces} \label{proofdensesob}
For the proof of Theorem \ref{densiteW}, we need two lemmas. Below,
for $I$ an open subset of $\T$, the symbol $\langle \cdot,\cdot\rangle_I$ 
stands for the duality bracket between 
$W_{\RR}^{-1/p,p}(I)$ and $W^{1-1/q,q}_{0,\RR}(I)$. We drop the subscript
when $I=\T$.

\begin{lemmasubsubsect} \label{selfadj}
For $\varphi\in W^{1-1/p,p}_{\R}(\T)$ and $\psi\in W^{1-1/q,q}_{\R}(\T)$, one has
\begin{equation}
\label{thetaconj}
\langle \partial_\theta (\mathcal{H}_\nu\varphi),\psi\rangle=\langle \varphi, 
\partial_\theta (\mathcal{H}_\nu\psi)\rangle.
\end{equation}
\end{lemmasubsubsect}
{\it Proof.} Assume first that $\varphi$ and $\psi$ lie in $C^{\infty}(\T)$. 
Let $f$ and $g$ be the solutions to \eqref{dbar}
on $\D$, normalized as in \eqref{normal:f}, such that 
$\mbox{Re }(\mbox{tr }f)=\varphi$ and $\mbox{Re }(\mbox{tr }g)=\psi$ respectively.
By definition of $\mathcal{ H}_\nu$, it holds that $\mbox{Im }(\mbox{tr
}f)=\mathcal{H}_\nu\varphi$ and $\mbox{Im
}(\mbox{tr }g)=\mathcal{H}_\nu\psi$. As pointed out in the proof of Proposition \ref{hnusob}, 
the functions $f$ and $g$ lie in
$W^{2,r}(\D)$  for all $r\in (1,\infty)$ by \cite[Thm 9.15]{gt}.
Put $u=\mbox{Re }f$ and $u_1=\mbox{Re }g$. By the divergence  formula
\[
0=\int_{\D} \left(\mbox{div}(\sigma \nabla u)u_1-\mbox{div}(\sigma \nabla
  u_1)u\right){dm}=\int_{\T} \sigma \, \left(
  (\partial_nu)u_1-(\partial_nu_1)u\right)\, { 
d\theta},
\]
and from \eqref{CRtang} we see that
$\sigma \partial_nu=\partial_\theta \mathcal{H}_\nu{\varphi}$ while
$\sigma \partial_nu_1=\partial_\theta \mathcal{H}_\nu{\psi}$.
This yields the desired conclusion for smooth $\varphi$ and $\psi$.
In the general case, pick two sequences  $(\varphi_k)_{k\geq 1}$,
$(\psi_k)_{k\geq 1}$ of smooth functions
converging respectively to $\varphi$
in $W^{1-1/p,p}(\T)$ and  to $\psi$ in $W^{1-1/q,q}(\T)$. 
By Corollary \ref{MRieszgen},
we have that $\mathcal{H}_\nu\varphi_k\rightarrow \mathcal{H}_\nu \varphi$ in $W^{1-1/p,p}(\T)$
and $\mathcal{H}_\nu\psi_k\rightarrow \mathcal{H}_\nu\psi$ in $W^{1-1/q,q}(\T)$,
therefore $\partial_\theta(\mathcal{H}_\nu\varphi_k)\rightarrow
\partial_\theta(\mathcal{H}_nu\varphi)$ in $W^{-1/p,p}(\T)$ and
$\partial_\theta(\mathcal{H}_\nu\psi_k)\rightarrow \partial_\theta(\mathcal{H}_\nu\psi)$ in
$W^{-1/q,q}(\T)$. Identity \eqref{thetaconj} now follows from the first part of the proof by a 
limiting argument. \hfill\fin\\

\noindent To proceed with the second lemma, we introduce the following
piece of notation that will be of use in the next section as well:
whenever $I \subset\T$ and  $J=\T \setminus I$ is the complementary subset, 
then for
$u_I$ (resp. $u_J$) a function on $I$ (resp. $J$) we put
$u_I\vee u_J$ for the concatenated
function on $\T$ which is $u_I$ (resp. $u_J$) on $I$
(resp. $J$).

Let now $I$, $J$ be proper open subsets of $\T$ such that
$J=\T\setminus\overline{I}\neq\emptyset$.
For any function $u_J\in
W^{1-1/p,p}_{0,\R}(J)$, we form the concatenated function $0\vee u_J$ and we
set 
\begin{equation}
\label{defAb}
A \, u_J= \partial_\theta \bigl(\mathcal{H}_\nu (0\vee
      u_J)\bigr)_{\vert_I} \, .
\end{equation}
Note that $0\vee u_J\in W^{1-1/p,p}_\R(\T)$, so that 
$A \, : \, W^{1-1/p,p}_{0,\R}(J) \to
W^{-1/p,p}_{\R}(I)$ 
is well-defined
and bounded by Corollary \ref{MRieszgen} and the boundedness of 
$\partial_\theta$ from $W^{1-1/p,p}(I)$ into $W^{-1/p,p}(I)$.
\begin{lemmasubsubsect}
\label{lem:Adenserange}
The operator $A$ defined in \eqref{defAb} has dense range. 
\end{lemmasubsubsect}
{\it Proof.} It is equivalent to show that the adjoint operator
$A^{\ast}:W^{1-1/q,q}_{0,\R}(I)\rightarrow W^{-1/q,q}_{\R}(J)$ is one-to-one. 
Now, for $u_I\in W^{1-1/q,q}_{0,\R}(I)$ and $u_J\in W^{1-1/p,p}_{0,\R}(J)$,
we get by Lemma \ref{selfadj} 

\[
\displaystyle \langle A^{\ast}u_I,u_J\rangle_J  =  \displaystyle \langle u_I,Au_J\rangle_I
=  \displaystyle  \langle u_I,\partial_\theta \bigl(\mathcal{H}_\nu
(0\vee u_J)\bigr)_{\vert_I} \rangle_I 
= \displaystyle \langle (u_I\vee 0),\partial_\theta \mathcal{H}_\nu
(0\vee
    u_J)\rangle 
\]
\[
=  \displaystyle \langle (0\vee u_J) ,\partial_\theta \mathcal{H}_\nu
(u_I\vee 0) \rangle 
=  \displaystyle  \langle u_J,\partial_\theta \Bigl(\mathcal{H}_\nu
(u_I\vee 0)\Bigr)_{\vert_J} \rangle_J,
\]
hence $A^{\ast}u_I=\partial_\theta \bigl(\mathcal{H}_\nu(u_I\vee
      0)\bigr)_{\vert_J}$.
Thus the relation  $A^{\ast}u_I=0$ means  that 
$\mathcal{H}_\nu(u_I\vee 0)$ is constant on each component of $J$.
Let $J_0$ be such a component and $f$
be the solution of \eqref{dbar} in $W^{1,p}(\D)$ such that 
$\mbox{Re }(\mbox{tr }f)=u_I\vee 0$, normalized as in \eqref{normal:f}
so that $\mbox{Im }(\mbox{tr }f)=\mathcal{H}_\nu(u_I\vee 0)$. Since
$u_I\vee0$ vanishes on $J$, there exists $c\in \R$ such that 
$\mbox{tr }f=ic$ on $J_0$. Therefore $f\equiv ic$ in $\D$ by
Proposition \ref{hardysobol} and Proposition \ref{cor:propHpnu} point (c).
In particular $\mbox{Re }(\mbox{tr }f)=0$ on $\T$, hence $u_I=0$ as desired.
\hfill\fin

\medskip

{\it Proof of Theorem \ref{densiteW}.} By definition of an extension set, 
we can write $I=\cup_j(a_j,b_j)$ where each $a_j$ lies at positive distance
from the $b_k$, therefore also from the $a_k$ for $k\neq j$.
This implies there is an open arc $I_1$ with $\overline{I_1}\neq\TT$
such that $\overline{I}\subset I_1$.
By the extension property and a smooth partition of unity argument,
each function in $W^{1-1/p}(I)$ extends to a function in $W_0^{1-1/p}(I_1)$.
Thus, upon trading $I$ for $I_1$, it is enough to prove that any function 
in $W^{1-1/p,p}_0(I)$ can be approximated in $W^{1-1/p,p}(I)$ by the trace 
of a solution of (\ref{dbar}). 

\noindent Let $\varepsilon>0$ and $\varphi_I\in W^{1-1/p,p}_0(I)$ 
with $\varphi_I=u_I+iv_I$,
where $u_I$ and $v_I$ are real valued. Set 
$v=\left(v_I-{\mathcal H}_{\nu}(u_I\vee 0)\right)_{\vert I}
\in W^{1-1/p,p}_{\R}(I)$.
By Lemma \ref{lem:Adenserange}, there exists $u_J\in W^{1-1/p,p}_{0,\R}(J)$ 
satisfying
\begin{equation} \label{approx4}
\left\Vert \partial_{\theta}v-\partial_{\theta}{\mathcal H}_{\nu}(0\vee u_J)\right\Vert_{W^{-1/p,p}(I)}\leq \varepsilon,
\end{equation}
from which it follows by elementary integration and the Poincar\'e inequality
on $I$  that, for some $c_v\in \R$ and some $C_I$ independent of $v$, 
\begin{equation} \label{approx3}
\left\Vert v-{\mathcal H}_{\nu}(0\vee u_J)-c_v\right\Vert_{W^{1-1/p,p}(I)}\leq C_I\varepsilon.
\end{equation}
Consider now the concatenated function $(u_I\vee u_J)\in W^{1-1/p}(\T)$, 
and define
$$
\psi=(u_I\vee u_J)+i{\mathcal H}_{\nu}(u_I\vee u_J)+ic_v.
$$
Then $\psi\in W^{1-1/p,p}(\T)$ by Corollary \ref{MRieszgen}, and
by construction it is the trace on $\T$ of a solution to (\ref{dbar}). Since
$$
\psi=(u_I\vee 0)+(0\vee u_J)+i{\mathcal H}_{\nu}(u_I\vee 0)+i{\mathcal H}_{\nu}(0\vee u_J)+ic_v,
$$
it follows from (\ref{approx3}) that
$$
\left\Vert \varphi_I-\psi\right\Vert_{W^{1-1/p,p}(I)}\leq C_I\varepsilon.
$$
\hfill\fin

\medskip

\subsubsection{Density for Hardy traces} \label{proofdensehard}
The proof of Theorem \ref{thm:densite} resembles that of 
Theorem \ref{densiteW}, but makes use of different operators.
Herafter, the symbol $\langle \cdot,\cdot\rangle$ 
indicates not only the duality between 
$W_{\RR}^{-1/p,p}(\T)$ and $W^{1-1/q,q}_{0,\RR}(\T)$, as in
Lemma \ref{selfadj}, but also the $L_\RR^p-L_\RR^q$ and the $W_\RR^{1,p}-W_\RR^{-1,q}$ duality.
This should cause no confusion.

First, we will need the following version of Lemma \ref{selfadj}:
\begin{lemmasubsubsect} \label{selfadjbis}
For all $\varphi\in W^{1,p}_{\R}(\T)$ and all $\psi\in L^q_{\R}(\T)$, one has
\[
\langle \partial_\theta ({\cal H}_\nu{\varphi}),\psi\rangle=
\langle \varphi, \partial_\theta ({\cal H}_\nu{\psi})\rangle.
\]
\end{lemmasubsubsect}
{\it Proof.} 
When $\varphi$ and $\psi$ are smooth the result is contained in
Lemma \ref{selfadj}, and the general case then follows  
by a limiting argument, using the continuiy properties of $\partial_\theta$ 
together with  Corollary \ref{MRieszgen} and Proposition \ref{hnusob}.
\hfill\fin\par
\medskip
From Lemma \ref{selfadjbis}, we deduce a relation between 
${\cal H}_\nu: L^p_{\RR}(\T)\to L^p_{\RR}(\T)$ and
its adjoint ${\mathcal H}_{\nu}^*:L^q_{\RR}(\T)\to L^q_{\RR}(\T)$.
\begin{lemmasubsubsect} \label{hnuadj}
For any function $u\in L^q_{\R}(\T)$, it holds that
$$
{\mathcal H}_{\nu}^{\ast}u=-\partial_{\theta}{\mathcal H}_{\nu}U,
$$
where $U\in W^{1,q}_{\R}(\T)$ is any function such that
\begin{equation} \label{defU}
\partial_{\theta}U=u-\frac 1{2\pi}\int_{0}^{2\pi} u(e^{i\theta}) d\theta.
\end{equation}
\end{lemmasubsubsect}
{\it Proof.} Let $u\in L^q_{\R}(\T)$, $v\in L^p_{\R}(\T)$, 
and consider $U\in W^{1,q}_{\R}(\T)$ such that (\ref{defU}) holds. 
Then, since ${\mathcal H}_{\nu}v$ has zero mean on $\T$,
we get from Lemma \ref{selfadjbis}
\[
\displaystyle \langle {\mathcal H}_{\nu}^{\ast}u,v\rangle =  \displaystyle \langle u,{\mathcal H}_{\nu}v\rangle
=   \displaystyle \langle u-\frac 1{2\pi}\int_{0}^{2\pi} u(e^{i\theta}) d\theta\,,\,{\mathcal H}_{\nu}v\rangle
\]
\[
 = \displaystyle \langle \partial_{\theta}U,{\mathcal H}_{\nu}v\rangle
 =  \displaystyle -\langle U,\partial_{\theta} {\mathcal H}_{\nu}v\rangle
 =  \displaystyle -\langle \partial_{\theta}{\mathcal H}_{\nu}U,v\rangle.
\]
\hfill\fin\par
\medskip
\begin{remarquesubsubsect}
If we restrict ${\mathcal H}_{\nu}$ to the space
of $L^q_\RR(\T)$-functions with zero mean, 
which is mapped into itself  by ${\cal H}_\nu$, we may summarize 
the content of Lemma \ref{hnuadj} as 
${\cal H}_\nu^*=-\partial_\theta{\cal H}_\nu\partial_{\theta}^{-1}$.
\end{remarquesubsubsect}
\noindent Let $I\subset \T$ be as in Theorem \ref{thm:densite} 
and put $J=\T\setminus I$.  For $u_J\in L^p_{\R}(J)$, let us define
\[
\cB \, u_J= \left. \bigl({\cal H}_\nu(0\vee
      u_J)\bigr)\right\vert_I. 
 \]
\begin{corollarysubsubsect}
\label{corthmLB}
The operator $\cB$ is bounded
from $L^p_{\R}(J)$ to $L^p_{\R}(I)/\R$ and has dense range.
\end{corollarysubsubsect}
Note that we consider $\cB$ as a mapping from
$L^p_\R(J)$ into the quotient space $L^p_\R(I)/\R$ rather than $L^p_\R(I)$.
That ${\mathcal B}$ has dense range means: given $v\in L^p_\R(I)$ 
and $\varepsilon>0$, there exist $u_J\in L^p_\R(J)$ and $c_v\in \R$ such that
\begin{equation} \label{approx}
\left\Vert v-{\mathcal H}_{\nu}(0\vee u_J)-c_v\right\Vert_{L^p(I)}\leq 
\varepsilon.
\end{equation}
\medskip

\noindent {\it  Proof.} 
Clearly $\cB$ is well-defined
and bounded by  Corollary \ref{MRieszgen}.
To prove it has dense range, it is enough to 
check that ${\mathcal B}^{\ast}:L^{q,0}_{\R}(I)\rightarrow L^q(J)$ 
is one-to-one; here, by $L^{q,0}_{\R}(I)$, we mean the subspace of $L^q_\R(I)$ 
comprised of functions with zero mean. 
For $\varphi_I\in L^{q,0}_\R(I)$, we deduce from
Lemma \ref{hnuadj} that
\begin{equation}
\label{derHPsi}
{\mathcal B}^{\ast}\varphi_I=-\left(\partial_{\theta}{\mathcal H}_{\nu}
\Psi\right)_{\vert J},
\end{equation}
where $\Psi\in W^{1,q}_{\R}(\T)$ is such that
\begin{equation}
\label{derPsi}
\partial_{\theta}\Psi=\varphi_I\vee 0.
\end{equation}
Consider the solution $f\in W^{1,p}(\D)$ to 
\eqref{dbar} satisfying
${\rm Re} \, (\tr\, f)= \Psi$ on $\T$, normalized
so that (\ref{normal:f}) holds. 
By \eqref{derPsi} it holds that $\partial_\theta\Psi=0$ a.e. on $J$,
and if  $\cB^{\ast}\varphi_I =0$ we get in addition from \eqref{derHPsi}
that $\partial_\theta {\mathcal H}_{\nu}\Psi= 0$ a.e.
on $J$. 
As $\tr\, f = \Psi+ i {\cal H}_\nu\Psi$,
it entails altogether that $\partial_\theta (\tr\,f)=0$ a.e. on $J$.
But from Theorem \ref{thmLB} point 
$(c)$ we know that $\partial f$ admits a non tangential 
limit a.e. on $\T$, and by \eqref{relim} we now see that 
$\tr\, \partial f =0$ a.e. on $J$. But 
Theorem \ref{thmLB} point $(b)$ and Lemma \ref{Beltcd} 
equation \eqref{BNRt} together imply that $\partial f=e^sF$, where $s$
is continuous on $\overline{\D}$ while $F\in H^p$. Necessarily then,
$\tr F=0$ on $J$ which has positive measure, hence $F\equiv0$ implying that
$\partial f=\overline{\partial} f\equiv0$. 
Therefore $f$ is  constant, in particular 
$0=\partial_\theta {\rm Re}(\tr f)=\partial_{\theta} \Psi
= \varphi_I \vee 0$ on $\T$. Therefore
$\varphi_I \equiv 0$ thus $\cB^*$ is injective, as desired.
\hfill\fin

\medskip

{\it Proof of Theorem \ref{thm:densite}.} 
Let $\varepsilon>0$ and $\varphi_I\in L^p(I)$ with $\varphi_I=u_I+iv_I$
where $u_I,v_I\in L^p_{\R}(I)$.
Set $v=v_I-{\mathcal H}_{\nu}(u_I\vee 0)$ and observe from (\ref{approx}) 
that
\begin{equation} \label{approx2}
\left\Vert v-{\mathcal H}_{\nu}(0\vee u_J)-c_v\right\Vert_{L^p(I)}\leq \varepsilon.
\end{equation}
for some $u_J\in L^p_\R(J)$ and $c_v\in \R$. Define
$$
\psi=(u_I\vee u_J)+i{\mathcal H}_{\nu}(u_I\vee u_J)+ic_v.
$$
By construction $\psi$ is  the trace on $\T$ of a $H^p_\nu$-function, and
since
$$
\psi=(u_I\vee 0)+(0\vee u_J)+i{\mathcal H}_{\nu}(u_I\vee 0)+i{\mathcal H}_{\nu}(0\vee u_J)+ic_v,
$$
it follows from (\ref{approx2}) that
$$
\left\Vert \varphi_I-\psi\right\Vert_{L^p(I)}\leq \varepsilon.
$$
The proof is now complete. \hfill\fin

\medskip
\subsection{A characterization of $(\tr \, H^{p}_\nu)^\perp$}
\label{sec:charac}
The expression for $\mathcal{H}_\nu^*$ obtained
in Lemma \ref{hnuadj} enables us to characterize the orthogonal space 
$(\tr \, H^{p}_\nu)^\perp$ of $\tr \, H^{p}_\nu$ for the duality product
\eqref{dualpair}, hereafter denoted by $\langle .,.\rangle$.

 {\it Proof of Proposition \ref{orthogformula}.} 
Let $\varphi = \varphi_1 + i \varphi_2$, with $\varphi_k \in L_\R^q(\T)$. 
In view of  Theorem \ref{thm:dirichletf}, 
 we get that
$\varphi \in (\tr \, H^{p}_\nu)^\perp$ if, and only if 
\begin{equation}
\label{orthogen}
\mbox{Re} \, 
\langle \varphi_1 + i \varphi_2 , u + i \cH_\nu u + i c \rangle = 0 \, , 
\ \  u \in  L_\R^p(\T),\  c \in \R \, .
\end{equation}
Picking for $u$ an arbitrary constant shows that
$\varphi_1$ and $\varphi_2$ have zero mean on $\T$, hence there are
$\Phi_1,\Phi_2\in W^{1,q}_\R(\T)$ such that $\Phi:=\Phi_1+i\Phi_2$
satisfies $\partial_\theta \Phi=\varphi$; we may impose in addition
that $\Phi$ itself has zero mean.
Now, from \eqref{orthogen} we get for every  $ u \in  L_\R^p(\T)$ that
\[
 0 = \mbox{Re} \, 
\langle u + i \cH_\nu u, \varphi_1 + i \varphi_2 \rangle 
= \langle u , \varphi_1 \rangle - \langle u , \cH_\nu^\ast  \varphi_2 \rangle, 
 \]
which means
$\varphi_1 = \cH_\nu^\ast  \varphi_2$ or else 
$\Phi_1 = - \cH_\nu \Phi_2$, by elementary integration and
using Lemma \ref{hnuadj}.
Therefore 
$\Phi = i \, (\Phi_2 + i \, \cH_\nu \Phi_2)$ lies in 
$i\tr\, H^q_\nu=\tr\, H^q_{-\nu}$, and since $\partial_\theta\Phi=\varphi$ 
we conclude that
\[
(\tr \, H^{p}_\nu)^\perp  \subset \partial_\theta \, \left(\tr\, H^{q}_{-\nu} \cap W^{1,q}(\T) \right) \, .
\]
Conversely let $\varphi= \partial_\theta \Phi$ for some
$\Phi\in \left(\tr\, H^{q}_{-\nu} \cap W^{1,q}(\T) \right)$.
We can write 
$\Phi = \gamma + i \, \cH_{-\nu} \, \gamma$ for some 
$\gamma\in W^{1,q}_\R(\T)$, and from Lemma \ref{hnuadj} applied with $-\nu$
rather than $\nu$ we obtain for $u \in  L_\R^p(\T)$ 
\[
\mbox{Re} \, \langle \varphi, u + i \cH_\nu u \rangle = 
\langle \partial_\theta \, \gamma , u \rangle - \langle \partial_\theta  \, (\cH_{-\nu} \, \gamma) , \cH_{\nu} \, u \rangle 
\]
\[
= 
\langle \partial_\theta \, \gamma , u \rangle + \langle \cH_{-\nu}^* (\partial_\theta  \, \gamma) , \cH_{\nu} \, u \rangle = \langle \partial_\theta \, \gamma , u \rangle + \langle \partial_\theta  \, \gamma , \cH_{-\nu} (\cH_{\nu} \, u )\rangle 
\]
\[= \langle \partial_\theta  \, \gamma , \int_0^{2 \pi} u d \theta \rangle = 0 \, ,
\] 
since $\cH_{-\nu} (\cH_{\nu} \, u ) = - u + \int_0^{2 \pi} u d \theta$
as follows immediately from the fact that $iH^p_\nu=H^p_{-\nu}$.
\hfill\fin

\medskip

\begin{remarquesubsect}
In order to extend the definition of $\cH_\nu$ to complex valued 
functions, it is natural to set ({\it cf.} \cite{ap, ap2})
\[
\cH_\nu (i \, u) = i \, \cH_{-\nu} (u) \, ,
\]
for we then have
\[
i \, u + i \, \cH_\nu (i \, u)  = i \, (u + i \, \cH_{-\nu} (u)) \, ,
\]
which is indeed a solution to \eqref{dbar}.
With this definition, we recap Proposition \ref{orthogformula} by saying 
that $\varphi \in  L^q(\T)$ belongs to $(\tr \, H^{p}_\nu)^\perp$ 
if and only if $\varphi = \partial_\theta \Phi$, where $\Phi\in W^{1,q}(\T)$ satisfies
\begin{equation}
\label{conjprop}
\cH_{-\nu} \Phi = - i \, \Phi , . 
\end{equation}
Note that if $\nu=0$ then  \eqref{conjprop} characterizes traces of holomorphic functions with zero mean. In this case \eqref{orthogen} 
follows readily from the Cauchy formula.
\end{remarquesubsect}

\section{Hardy spaces over Dini-smooth domains}
\label{dini-sec}
In this final section, we indicate how the spaces $H_\nu^p$ that we
studied on the disk can be defined 
more generally on Dini-smooth domains.\par

\medskip

\noindent \Rd Recall that a function $h$ is called {\emph {Dini-continuous}} if 
$\int_0^\varepsilon(\omega_h(t)/t)dt<+\infty$ for some, 
hence any $\varepsilon>0$, 
where $\omega_h$ is the modulus of continuity of $h$. 
A function is said to be {\emph {Dini-smooth}} if it has Dini-continuous 
derivative.
A bounded planar domain in $\C$ is termed Dini-smooth if its 
boundary is a Jordan 
curve with nonsingular Dini-smooth parametrization. Such domains $\Omega$
have the  property that any conformal map from $\D$ onto $\Omega$ extends 
continuously from  $\overline{\D}$  onto $\overline{\Omega}$ together with its 
derivative, in such a way that the latter is never zero 
\cite[thm 3.5]{Pommerenke}, \Rd and that is why we are able to generalize our previous results to this class of domains. \Bk \par

\medskip

\noindent Let $\Omega \subset \C$ be a simply connected Dini-smooth domain, as defined in the introduction, and assume that $\nu \in W^{1, \infty}(\Omega)$ with 
$\|\nu\|_{L^\infty(\Omega)} \leq \kappa <1$.
Any conformal transformation $\psi$ from $\D$ onto $\Omega$ has a 
$C^1(\overline{\D})$ conformal extension onto $\overline{\Omega}$
that we still denote by $\psi$ \cite[thm 3.5]{Pommerenke}.
Introduce the Hardy classes $H^p_\nu(\Omega)$ as the space of functions 
$f$ on $\Omega$ such that $f \circ \psi \in H^p_{\nu \circ \psi}(\D)$.
A straightforward computation \cite[Ch. 1, C]{ahlfors}
shows that this class does not depend on the 
particular choice of $\psi$  and consists of distributional solutions to
\eqref{dbar}. Note that  
$\nu \circ \psi \in W^{1, \infty}(\D)$ with $\|\nu\circ \psi \|_{L^\infty(\D)} \leq \kappa <1$. Similarly,  the 
class $G^p_\alpha(\Omega)$ consists of those functions $g$ on $\Omega$ 
such that $g \circ \psi \in G^p_{\alpha \circ \psi}(\D)$. 

As in the classical case of holomorphic Hardy spaces over 
simply connected domains \cite{duren}, it appears that 
$w \in G^p_\alpha(\Omega)$ if and only if it is a solution to  (\ref{eq:w}) 
such that $|w|^p$ has a harmonic majorant (we insist that harmonic is 
understood here in the usual sense). Indeed, it is enough to check this property in $\D$, since it is preserved by composition with conformal maps. Then, 
it is a consequence of Theorem \ref{trick2} that functions in 
$G^p_\alpha(\D)$ possess this property. For the converse, 
observe that any solution $w$ to (\ref{eq:w}) can be factorized as 
$w = e^s \, F$, where $s \in C(\overline{\D})$ and $F$ is holomorphic
in $\D$. Since $|F|^p = e^{- p \mbox{\small Re} \, s} \, |w|^p$, it admits
a harmonic majorant, and therefore it belongs to $H^p(\D)$, which
immediately yields that $w \in G^p_\alpha(\D)$.  

Another possibility to generalize $H^p_\nu$ to a Dini-smooth
simply connected domain $\Omega$ is to parallel the definition
of the so-called Smirnov classes \cite[Ch. 10]{duren}.
Namely, one requires the uniform boundedness of the $L^p$ norm 
on a sequence of rectifiable contours eventually encompassing
every compact subset of $\Omega$. For Dini-smooth domains,
the two generalizations turn out to be equivalent.

The results of Sections \ref{gpalphaprop}, \ref{kkkkk}, \ref{lllllllll}, and 
\ref{mmmmmmmmmmmmmm} remain valid in this framework,
as can be seen easily by composing with $\psi$ and appealing to the
regularity of $\psi^{-1}$.

\section{Conclusion}
\label{conclu}
This paper took a few steps towards a theory of
two-dimensional Hardy spaces for the conjugate Beltrami equation
on simply connected domains. Conspicuously missing is a factorization theory, 
whose starting point should be a characterization of those pairs $(s, F)$
for which \eqref{factBN} holds. Also, we did not pursue the solution to
the extremal problems
stated in Theorem \ref{thm:dualite}, points $(ii)-(iii)$. Finally, the case
of multiply connected domains, that motivated in part the present 
study ({\it cf.} the introduction), was not \Gr touched upon. \Bk It is to be hoped that
suitable deepenings of the present techniques will enable one to
approach such issues.

\appendix
\section{{\Bl Appendix.} \Rd Proof of technical results \Bk} \label{technic}
\Rd We establish here Proposition \ref{properties} and Lemma \ref{improve}. \par

\medskip

\noindent {\bf Proof of Proposition \ref{properties}: } \Bk The boundedness of $\mathcal C$ from $L^p(\T)$ onto
$H^p(\D)$ follows from the M. Riesz theorem, mentioned already. 
To establish the
second half of 1, let $\varphi \in L^p( \T)$. From 
the definition of the $H^p$-norm and the M. Riesz theorem, 
we get for each $r \in (0,1)$ and every $\epsilon > 0$ that
\[
\|{\cal C} \varphi \|_{L^p(\T_r)} \leq C_p \, \|\varphi\|_{L^p( \T)} \leq
    \frac{C_p}{(1-r)^{\epsilon}} \|\varphi\|_{L^p( \T)} \, ,
\]
where the last inequality is trivial. It then follows from \cite[thm
5.5]{duren} that
\begin{equation}
\label{boundC}
\|\partial{\cal C} \varphi \|_{L^p(\T_r)} \leq
\frac{C_{\epsilon,p}} {(1-r)^{1+\epsilon}} \|\varphi\|_{L^p( \T)} \, .
\end{equation}
If moreover $\varphi \in W^{1,p}( \T)$, integrating by parts the Cauchy 
formula gives us
\[\partial {\mathcal C} \varphi(z)=\frac{1}{2i\pi}\int_\T
\frac{\varphi(\xi)}{(\xi-z)^2}\,d\xi=
\frac{1}{i}\int_\T
\frac{\partial_t \varphi(\xi)}{\xi(\xi-z)}\,d\xi,
\]
from which it follows by the M. Riesz theorem again that
\begin{equation}
\label{boundC2}
\|\partial {\cal C} \varphi \|_{L^p(\T_r)} \leq C_p
\|\varphi\|_{W^{1,p}(\T)} \, .
\end{equation}
Introduce the operators $A_r^1:
W^{1,p}( \T) \to L^p(\TT_r)$ and  $A_r^0: L^{p}( \T) \to L^p(\TT_r)$, 
where in both cases $A_r^i \varphi = (\partial {\cal C} \varphi)_{|_{\T_r}}$.
We gather from (\ref{boundC2}) and (\ref{boundC}) that
\[
|\!|\!|A_r^1 |\!|\!| \leq C_p \ \mbox{ while }\  |\!|\!|A_r^0 |\!|\!| \leq
\frac{C_{\epsilon,p}}{(1-r)^{1+\epsilon}} \, .
\]
Since $W^{1-1/p,p}( \T)=[W^{1,p}(\T),L^p(\T)]_{1/p}$, 
interpolating between $A^1_r$ and $A^0_r$ yields
\[
\|\partial {\cal C} \varphi \|_{L^p(\T_r)} \leq
\frac{C_{\epsilon,p}^\prime}{(1-r)^{(1+\epsilon)/p}}
\|\varphi\|_{W^{1-1/p,p}( \T)} \, .
\]
Choosing $\epsilon$ so small that $p> 1+\epsilon$, the above right hand side 
is integrable w.r.t. $r \in
(0,1)$, so by Fubini's theorem $\partial {\cal C}
\varphi \in L^p(\D)$ as soon as $\varphi \in W^{1-1/p,p}( \T)$. 
Since $\bar\partial {\cal C} \varphi = 0$, this establishes 1.

\par
\noindent Assertion 2 
is a consequence of the fact that $S$ is a $L^2(\C)$-isometric
Calder\'on-Zygmund operator, see \cite[Ch. V, Sec. D]{ahlfors} or
\cite[Ch. II, { Thm} 3]{stein}.
\par
\noindent Next, observe that if $K\subset\C$ 
and $w \in L^p(\C)$, we have for $z\in K$ 
$$ \breve Tw(z)=\frac1{\pi}\left(
(\chi_{\D}w)*g_K\right)(z),~~\mbox{with}~~g_K(\xi):={\chi_{K+\D}(\xi)\over \xi},
$$ where $\chi_E$ denotes the characteristic function of a set $E$.
Clearly $g_K\in L^1(\C)$ when $K$ is compact, showing
that $\breve{T}$ is bounded from $L^p(\C)$ into $L^p_{loc}(\C)$. 
We claim that $\bar\partial\breve{T}w=\chi_\D w$ and
$\partial\breve{T}w=S(\chi_\D w)$ in the sense of distributions.
When $w$ is $C^2$-smooth and compactly supported in $\D$, whence
$\chi_\D w=w$, this is 
a simple computation \cite[Ch. V, Lem. 2]{ahlfors}. 
In the general case, pick a sequence of functions
$v_n\in \cD(\D)$ converging to
$\chi_\D w$ in $L^p(\C)$. By what precedes $\breve{T}v_n$ converges to
$\breve{T}(\chi_\D w)=\breve{T}w$ in $L^p_{loc}(\C)$, hence as distributions
\begin{equation}
\label{limdbarw}
\bar\partial(\breve{T}w)=\lim_n\bar\partial(\breve{T}v_n)=\lim_nv_n=
\chi_\D w,
\end{equation}
where the last limit holds in $L^p(\C)$ thus {\it a fortiori} in the 
distributional sense.
Using the $L^p(\C)$ boundedness of $S$, a similar argument yields
\begin{equation}
\label{limdw}
\partial (\breve{T}w)=\lim_n\partial(\breve{T}v_n)=\lim_nS(v_n)=
S(\chi_\D w),
\end{equation}
proving the claim. Since  $\left\Vert
\partial \breve{T}w\right\Vert_{L^p(\C)}=\left\Vert
S(\chi_\D w)\right\Vert_{L^p(\C)}
\leq C\left\Vert
w\right\Vert_{L^p(\D)}$, assertion $3$ holds.
\par
\noindent As to assertion $4$, the boundedness of $T$ from $L^p(\D)$
to $W^{1,p}(\D)$ follows from $3$ and the identity 
$(\breve{T} \breve{w})_{|_\D} = T w$, an so do the relations
$\bar\partial T w =w$ and $\partial T w =(S\breve{w})_{|_\D}$
in view of \eqref{limdbarw} and \eqref{limdw}.
The compactness of $T:L^p(\D)\rightarrow L^p(\D)$
then follows from the compactness of the injection
$W^{1,p}(\D)\rightarrow L^p(\D)$ \cite[Ch. VI, { Thm} 6.2]{Adams}. 

Finally, set $T_{\alpha}w=T(\alpha \overline{w})$. 
By a theorem of F. Riesz,
to see that $I-T_{\alpha}$ is an isomorphism of $L^p(\D)$, it is
enough since $T_{\alpha}$ is compact to check that 
$I-T_{\alpha}$ is one-to-one. Let $w\in
L^p(\D)$ be such that $(I-T_{\alpha})w=0$ and
set $u=\breve{T}(\breve{\alpha \overline{w}}) \in W^{1,p}_{loc}(\C)$. 
Observe there is $p_1>2$ such that 
$u\in L^{p_1}_{loc}(\C)$. Indeed, by the Sobolev imbedding theorem,
if $p<2$ then  $W^{1,p}_{loc}(\C)$ is embedded in $L^{p^{\ast}}_{loc}(\C)$,
with $p^{\ast}=\frac{2p}{2-p}$, whereas if $p\geq2$ then 
$W^{1,p}_{loc}(\C)$ is
embedded in $L^\lambda_{loc}(\C)$ for every $\lambda\in (2,\infty)$. 
\par
\noindent Now, since $w=T_{\alpha}w$,
we have $u=w$ in $\D$ and so $w\in L^{p_1}(\D)$ hence
$\breve{\alpha \overline{w}}\in L^{p_1}(\C)$.
Thus, by assertion $3$, we get that $u\in
W^{1,p_1}_{loc}(\C)$. Moreover, from \eqref{limdbarw} and 
since $u=w$ in $\D$,
it holds in the sense of distributions that
\begin{equation} \label{dbaru1}
\bar\partial u=\breve{\alpha 
\overline{w}}=\breve{\alpha} \, \overline{u}~~\mbox{ a.e.  in
}\C.
\end{equation}
In addition, 
$u(z)$ clearly goes to $0$ when $\left\vert
z\right\vert$ goes to $+\infty$. It now follows from the generalized 
Liouville theorem \cite[Prop. 3.3]{ap} 
that $u=0$, therefore $w=0$.
 \hfill\fin\par
\medskip

\Rd {\it Proof of Lemma \ref{improve}: }  \Bk It is a result by
Hardy and Littlewood \cite[thm 5.9]{duren} that for any $g\in H^p$
\begin{equation}
\label{HLf}
\|g\|_{L^{p_1}(\T_r)} \leq C_{p} \|g\|_{L^{p}(\T_r)}
(1-r)^{1/p_1-1/p} \, , \ 0 < r < 1 \, ,\ \ p\leq p_1\leq\infty\,.
\end{equation}
Taking $p_1\in[p,2p)$ and raising \eqref{HLf} to the power $p_1$,
we obtain upon integrating with respect to $r$
\[
\|g\|_{L^{p_1}(\D)}^{p_1} 
\leq
C_{p,p_1} \|\tr\, g\|^{p_1}_{L^p(\T)}  
\]
thereby proving $1.$ For the rest of the proof, we fix $p_1\in(2,2p)$.

\noindent Assume now that $w - T(\alpha \overline{w}) = g \in H^p$. 
By Proposition \ref{properties}, point 4,
we get $T(\alpha \overline{w}) \in W^{1,p}(\D)$. Moreover, as 
already stressed in the proof of that proposition, we know from the Sobolev 
imbedding theorem  that $W^{1,p}(\D)\subset L^{p_2}(\D)$ for
some $p_2 >2$. Thus, by point 1 and H\"older's inequality, 
$w = g + T(\alpha\overline{w})$ belongs to $L^{p^*}(\D)$ with 
$p^* = \min \{p_1, p_2\} >2$, and consequently
$T (\alpha\overline{w}) \in W^{1,p^*}(\D) \subset
C^{0,1-2/p^*}(\overline{\D})$ by the Sobolev imbedding theorem again.
Using Proposition \ref{properties}, point 3,
we establish similarly that 
$\breve{T}\left(\breve{\alpha \overline{w}}\right) \in W_{loc}^{1,p^*}(\C) 
\subset C_{loc}^{0,1-2/p^*}(\C)$. 
Recaping what we just did, we obtain
$$
\|T(\alpha\overline{w})\|_{W^{1,p^*}(\D)}\leq c_{p,\alpha}\|w\|_{L^{p^*}(\D)}=
c_{p,\alpha}\|g+T(\alpha\overline{w})\|_{L^{p^*}(\D)}
\leq c_{p,\alpha}'\left(\|T(\alpha\overline{w})\|_{L^{p_2}(\D)}+\|g\|_{L^{p_1}(\D)}\right)
$$
$$
\leq c_{p\alpha}''\left(\|T(\alpha\overline{w})\|_{W^{1,p}(\D)}+\|g\|_{H^p(\D)}\right)
\leq c_{p,\alpha}'''\left(\|w\|_{L^{p}(\D)}+\|g\|_{H^p(\D)}\right)
$$
which is (\ref{td}).
\hfill\fin\par

\bigskip

\noindent{\small{{\bf { Aknowledgments. } }}}

We wish to thank K. Astala, S. Hofmann, F. Murat and L. 
P\"aiv\"arinta for helpful 
discussions and advice. \Gr We are indebted to the referee whose remarks helped us to improve the presentation of this paper. \Bk This work was partially supported by the ANR project ``AHPI'' (ANR-07-BLAN-0247-01).


\begin{thebibliography}{AAA}
\bibitem{AAK} V. M. Adamjan, D. Z. Arov, M. G. Krein,
Analytic properties of {S}chmidt pairs for a {H}ankel
operator and the generalized {S}chur--{T}akagi problem,
{\it Math. USSR Sbornik}, { 15}, 31--73, 1971.
\bibitem{Adams} R. Adams, {\it Sobolev spaces}, Academic Press, 1975.
\bibitem{ahlfors} L. Ahlfors, {\it Lectures on quasiconformal
mappings}, Wadsworth and Brooks/Cole Advanced Books and Software,
Monterey, CA, 1987.
\bibitem{abl} D. Alpay, L. Baratchart, J. Leblond, Some extremal
problems linked with identification from partial frequency data, {\it
10th conference in Analysis and optimization of systems, Sophia
Antipolis, 1992}, in J. L. Lions, R. F. Curtain, A. Bensoussan,
Lecture notes in Control and Information Sc., { 185}, Springer
Verlag, 563-573, 1993.
\bibitem{ast} K. Astala, Area distortion of quasi-conformal mappings,
{\it Acta Math.} { 173}, 37-60, 1994.
\bibitem{aim} K. Astala, T. Iwaniec, G. Martin, {\it Elliptic Partial
Differential Equations and Quasiconformal mappings in the plane},  Princeton Math. Series, 2008.
\bibitem{ap} K. Astala, L. P\"aiv\"arinta, Calder\'on's inverse
conductivity problem in the plane, {\it Ann. of Math.} (2) { 16},
no. 1, 265--299, 2006.
\bibitem{ap2} K. Astala, L. P\"aiv\"arinta, A boundary integral equation for Calder\'on's inverse conductivity problem, 
{\it Proc. 7th Int. Conf. on Harmonic Analysis and
PDEs., Madrid (Spain), 2004}, Collect. Math., 127--139, 2006.
\bibitem{ABLP} B. Atfeh, L. Baratchart, J. Leblond, J. R. Partington,
Bounded extremal and Cauchy Laplace problems on the sphere and shell,
{\it submitted for publication}.
\bibitem{aq} P. Auscher, M. Qafsaoui, {\it Observations on $W^{1,p}$ estimates for divergence elliptic equations with VMO coefficients, } Boll. Unione 
Mat. Ital. Sez. B, Artic. Ric. Mat. { 8},  5  (2002),  no. 2, 487--509. 
\bibitem{BBHL} L. Baratchart, A. Ben~Abda, F. Ben~Hassen, J. Leblond,
Recovery of pointwise sources or small inclusions in 2D domains and 
rational approximation, {\it Inverse Problems}, { 21}, 51-74,
2005.
\bibitem{BLPprep} L. Baratchart, J. Grimm, J. Leblond,
J.R. Partington, Asymptotic estimates for interpolation 
and constrained
  approximation in {$H^2$} by diagonalization of 
Toeplitz operators, {\it Integral equations and operator theory},
{ 45}, 269--299, 2003.
\bibitem{bl} L. Baratchart, J. Leblond, Hardy approximation to $L^p$
functions on subsets of the circle with $1\leq p<\infty$, {\it
Constr. Approx.} { 14}, 41-56, 1998.
\bibitem{BLMS} L. Baratchart, J. Leblond, F. Mandr\'ea, E. B. Saff,
How can the meromorphic approximation help to solve some 
2{D} inverse problems for the {L}aplacian?,
{\it Inverse Problems}, { 15}, 79--90, 1999.
\bibitem{blp2} L. Baratchart, J. Leblond, J. R. Partington, Hardy
approximation to $L^{\infty}$ functions on subsets of the circle, {\it
Constr. Approx.} { 12}, 423-436, 1996.
\bibitem{blp3} L. Baratchart, J. Leblond, J. R. Partington, Problems
of Adamjan-Arov-Krein type on subsets of the circle and minimal norm
extensions, {\it Constr. Approx.},  {16}, 333-357, 2000.
\bibitem{blpt} 
L. Baratchart, J. Leblond, J.R. Partington, N. Torkhani,
Robust identification in the disc algebra from band-limited data,
{\it IEEE Trans. on Automatic Control}, { 42} (9), 1997.
\bibitem{BMSW06}
L. Baratchart, F. Mandr\`ea, E. B. Saff, F. Wielonsky,
2-{D} inverse problems for the {L}aplacian: a meromorphic approximation 
approach, {\it J. Math. Pures Appl.}, { 86}, 1-41, 2006.
\bibitem{BarSeyf} L. Baratchart,  F. Seyfert,
An $L^p$ analog to the {A}{A}{K} theory,
{\it Journal of Functional Analysis}, { 191}, 52-122, 2002.
\bibitem{bn} L. Bers, L. Nirenberg, On a representation theorem for
linear elliptic systems with discontinuous coefficients and its
applications, {\it Convegno internazionale sulle equazioni derivate e
parziali}, Cremonese, Roma, 111-138, 1954.
\bibitem{blum} J. Blum,
{\it Numerical Simulation and Optimal Control in Plasma Physics: With 
Applications to Tokamaks (Modern Applied Mathematics)}
John Wiley \& Sons, 1989.
\bibitem{boj} B. Bojarski, Homeomorphic solutions of Beltrami systems, {\it Dokl. Akad. Nauk. SSSR} { 102}, 661-664, 1955.
\bibitem{BW}
J. Bourgain, T. Wolff.
A remark on gradients of harmonic functions in dimension $d \geq 3$.
{\em Colloq. Math.}, 60/61, 253--260, 1990.
\bibitem{campanato} S. Campanato, Elliptic systems in divergence
form, Interior regularity, Quaderni, Scuola Normale Superiore Pisa, 1980.
\bibitem{carleson} L. Carleson, {\it Selected problems in exceptional sets}, Van Nostrand Math. Studies, vol. {\bf 13}, Princeton, 1967.
\bibitem{cp} I. Chalendar, J. R. Partington, Interpolation between
Hardy spaces on circular domains with applications to
approximation. {\it Arch. Math. (Basel)} { 78}, 223-232, 2002.
\bibitem{cps} I. Chalendar, J. R. Partington, M. Smith, Approximation
in reflexive Banach spaces and applications to the invariant subspace
problem, {\it Proc. A. M. S.} { 132}, 1133-1142, 2004.
\bibitem{cfmoz} A. Clop, D. Faraco, J. Mateu, J. Orobitg, X. Zhong,
Beltrami equations with coefficients in the Sobolev space $W^{1,p}$,
{\it Publ. Mat.}, 53 (1), 197-230, 2009.
\bibitem{dav} G. David, Solutions de l'\'equation de Beltrami avec
$\left\Vert \mu\right\Vert_{\infty}=1$, {\it
Ann. Acad. Sci. Fenn. Ser. A1 Math.} { 13}, 25-70, 1988.
\bibitem{Demengel} F. Demengel and  G. Demengel,
{\it Espaces fonctionnels pour la th\'eorie des \'equations},
Savoirs Actuels, Edp Sciences, 2007.
\bibitem{DFT} J. C. Doyle, B. A. Francis, A. R. Tannenbaum,
{\it Feedback Control Theory},
Macmillan Publishing Company, 1992.
\bibitem{duren} P. L. Duren, {\it Theory of $H^p$ spaces}, Pure and
Applied Mathematics, Vol. 38 Academic Press, New York-London, 1970.
\bibitem{FJR}
E. Fabes, M. Jodeit, N. Rivi\`ere,
Potential techniques for boundary value problems on {$C^1$} domains,
{\it Acta Math.}, { 141}, 165--186, 1978.
\bibitem{Forster} O. Forster, {\it Lectures on Riemann surfaces},
G.T.M.  81, Springer. 1981.
\bibitem{gar} J. Garnett, {\it Bounded analytic functions}, Pure and
Applied Math. { 96}, Academic Press, 1981.
\bibitem{gt} D. Gilbarg, N. Trudinger, {\it Elliptic partial
differential equations of second order}, Springer Verlag, 1983.
\bibitem{Glover}
K. Glover, All optimal Hankel--norm approximations of linear multivariable
                 systems and their $L^\infty$--error bounds,
{\it Int. J. Control }, { 39 (6)}, 1115--1193, 1984.
\bibitem{gk} G. M. Goluzin, N. I. Krylov, Generalized Carleman formula
and its applications to analytic extension of functions, {\it
Mat. Sb.} { 40}, 144-149, 1933.
\bibitem{Grisg} P. Grisvard,
{\it Elliptic problems in nonsmooth domains},
Pitman, 1985.
\bibitem{gmsv} V. Gutlianskii, O. Martio, T. Sugawa, M. Vuorinen, On
the degenerate Beltrami equations, {\it Trans. Amer. Math. Soc.} {
357}, 875-900, 2005.
\bibitem{Hormander}L. H\"ormander, 
{\it Complex analysis in several variables}, North Holland, 1990.
\bibitem{imsurv} T. Iwaniec, G. Martin, What's new for the Beltrami
equation? in {\it National Research Symposium on Geometric Analysis
and Applications}, Proc. of the CMA, { 39}, 2000.
\bibitem{im2} T. Iwaniec, G. Martin, {\it Geometric Function Theory
and Non-linear Analysis}, Oxford Univ. Press, 2001.
\bibitem{im1} T. Iwaniec, G. Martin, The Beltrami Equation, {\it Mem. A.M.S.}, vol. 191, 893, 2007. 
\bibitem{imh2} 
B. Jacob, J. Leblond, J.-P. Marmorat, J.R. Partington,
A constrained approximation problem arising in  
     parameter identification,
{\it Linear Algebra and its Applications},
{ 351-352}, 487--500, 2002.
\bibitem{LJMP}
M. Jaoua, J. Leblond, M. Mahjoub, J.R Partington,
Robust numerical algorithms based on analytic approximation for the solution 
of inverse problems in annular domains,
{\it IMA J. of Applied Math.},  to appear.
\bibitem{JeKe} D. Jerison, C. Kenig,
Boundary Value Problems on {L}ipschitz Domains,
in {\it Studies in Partial Differential Equations}, 1-68, MAA Stud. Math
{ 23}, AMS, 1982.
\bibitem{KN} M.G. Krein, P.Y. Nudel'man,
Approximation of $L^2(\omega_1, \omega_2)$ functions by minimum--
 energy transfer functions of linear systems,
{\it Problemy Peredachi Informatsii}, Eng. transl., { 11 (2)},
37--60, 1975.
\bibitem{lv} O. Lehto, K. I. Virtanen, {\it Quasiconformal mappings in
the plane}, Second Edition, Springer Verlag, New York, 1973.
\bibitem{LioMag1}
J. L. Lions, E. Magenes,
{\it Probl\`emes aux limites non homog\`enes et applications},
vol. 1, Dunod, 1968.
\bibitem{m} Y. Mizuta, Existence of various boundary limits of Beppo Levi functions of higher order, {\it Hiroshima Math. J.} {\bf 9}, 717-745, 1979.
\bibitem{mor} C. B. Morrey, On the solutions of quasi-linear elliptic
partial differential equations, {\it Trans. Amer. Math. Soc.} {
43}, 1, 126-166, 1938.
\bibitem{Nikolskiis}
N. K.  Nikolskii, {\it Treatise on the shift operator}, 
Grundl. der Math. Wissenschaften { 273}, Springer,
1986.
\bibitem{Nikolskii} N. K.  Nikolskii,
{\it Operators, functions, and systems: an easy reading}, 
Mathematical surveys and monographs, { 92-93}, Amer. Math. Soc.,
2002.
\bibitem{Parfenov} O. G. Parfenov, 
Estimates of the singular numbers of a {Carleson} operator, 
{\it Math USSR Sbornik}, { 59(2)}, 497-514, 1988.
\bibitem{patil} D. J. Patil, Representation of $H^p$ functions, {\it
Bull. Amer. Math. Soc.} { 78}, 4, 617-620, 1972.
\bibitem{Peller} V. V. Peller,
{\it Hankel Operators and their Applications}, 
Springer, 2003.
\bibitem{Pommerenke} Ch. Pommerenke, {\it Boundary behaviour of conformal maps}, Springer Verlag, 1992.
\bibitem{Prokhorov2002} V. A. Prokhorov,
On {$L^p$}-generalization of a theorem of {A}damjan, {A}rov, and {K}rein,
{\it J. Approximation Theory}, { 116(2)}, 380--396, 2002.
\bibitem{Ransford}
T. Ransford, {\it Potential theory in the complex plane},
Student Texts 28, London Math. Soc., 1995.
\bibitem{r} Yu. G. Resetnyak, Boundary behaviour of functions with generalized derivatives, {\it Siberian Math. J. } {\bf 13}, 285-290, 1972.
\bibitem{Schwartz} L. Schwartz,
{\it Th\'eorie des distributions}, Hermann, 1978.
\bibitem{seyfertIMS2003}
F. Seyfert, J. P. Marmorat, L. Baratchart, S. Bila, 
J. Sombrin,
Extraction of Coupling Parameters For Microwave Filters: 
Determination of a Stable Rational Model from Scattering Data,
{\it Proceedings of the International Microwave Symposium, 
       Philadelphia}, 2003.
\bibitem{stein} E. M. Stein, {\it Singular integrals and
differentiability properties of functions}, Princeton Univ. Press,
1970.
\bibitem{SteinWeiss} 
E.M. Stein, G. Weiss, 
{\it {Introduction to Fourier analysis on euclidean spaces}},
Princeton Univ. Press, 1971.
\bibitem{Torchinsky}
A. Torchinsky, {\it {Real-variable Methods in Harmonic Analysis}},
Academic Press, San Diego, 1986.
\bibitem{Ziemer}
W. P. Ziemer, {\it {Weakly Differentiable Functions}},
Graduate Texts in Mathematics, {\bf 120}, 
Springer-Verlag, New York, 1989.
\end{thebibliography}
\end{document}